%% file: main.tex
\documentclass[a4paper,fleqn]{cas-sc}

\usepackage{amsmath,amssymb,amsthm,mathtools,bm}
\usepackage[numbers,sort&compress]{natbib}
\usepackage{acro}
\usepackage{algorithm,algpseudocode}
\usepackage{bbm}
\usepackage{graphicx,booktabs,subcaption}
\usepackage{etoolbox}
\preto{\section}{\acresetall}
\usepackage[nameinlink,capitalise,noabbrev]{cleveref}
\usepackage[textsize=footnotesize,colorinlistoftodos]{todonotes}
\usepackage{xspace}

\input{acronyms}
\input{definitions}

\begin{document}

\shorttitle{Model-guided design of nanoparticle dispersity}
\shortauthors{Gilch et~al.}

\title[mode=title]{Model-Guided Design of Nanoparticle Dispersity through Process-Chain Optimization}

\author[1]{Andrea Gilch}[orcid=0009-0000-9267-4077]
\ead{andrea.gilch@fau.de}
\credit{Conceptualization, Methodology, Software, Formal analysis, Data curation, Investigation, Validation, Visualization, Writing -- original draft, Writing -- review \& editing}

\author[2]{Nabi E. Traoré}[orcid=0000-0003-1279-2503]
\ead{nabi.traore@fau.de}
\credit{Conceptualization, Methodology, Data curation, Writing -- original draft, Writing -- review \& editing}

\author[1]{Frauke Liers}[orcid=0000-0002-1671-0344]
\ead{frauke.liers@fau.de}
\credit{Conceptualization, Methodology, Project administration, Supervision, Writing -- review \& editing}

\author[3]{Lukas Pflug}[orcid=0000-0001-8001-5832]
\cormark[1]
\ead{lukas.pflug@fau.de}
\credit{Conceptualization, Formal analysis, Methodology, Software, Supervision, Writing -- review \& editing}

\cortext[1]{Corresponding author}

\affiliation[1]{
  organization={Department of Data Science and Department of Mathematics, Friedrich-Alexander-Universit\"at Erlangen-N\"urnberg},
  city={Erlangen},
  country={Germany}}

\affiliation[2]{
  organization={Institute of Interfaces and Particle Technology, Friedrich-Alexander-Universit\"at Erlangen-N\"urnberg},
  city={Erlangen},
  country={Germany}}

\affiliation[3]{
  organization={Competence Center Scientific Computing, Friedrich-Alexander-Universit\"at Erlangen-N\"urnberg},
  city={Erlangen},
  country={Germany}}

\begin{abstract}
The design of high-quality nanoparticles requires the production of materials with narrow particle size distributions, as these strongly influence the functionalities of the final product. To achieve this goal, it is useful to consider a coupled process chain in which nanoparticles are first synthesized and subsequently fractionated by chromatographic separation. We introduce a novel mathematical optimization model for the integrated optimization of coupled nanoparticle synthesis and separation processes. The proposed model simultaneously accounts for process performance and economic objectives through a cost function that incorporates both production costs and achievable profit defined by the particle quality. For efficient evaluation of the synthesis model, we employ the exact method of moments. The resulting non-convex process-design problem is solved numerically using local optimization. Optimization results of our approach are demonstrated for the prototypical design of gold nanoparticles. The computational results show that optimized solutions of the coupled synthesis--separation process are capable of accurately representing both process stages while achieving optimized particle quality. The proposed model highlights the potential of mathematical optimization for the design of coupled nanoparticle production process-chains.
\end{abstract}


\maketitle

\input{1_Introduction}
\input{2_Simulation}
\input{3_Optimization}
\input{4_Conclusion}

\section*{Acknowledgments}
The authors acknowledge funding  by the Deutsche Forschungsgemeinschaft (DFG, German Research Foundation), Project-ID 416229255, SFB 1411.

\printcredits

\section*{Declaration of competing interest}
The authors declare that they have no known competing financial interests
or personal relationships that could have appeared to influence the work
reported in this paper.

\section*{Declaration of generative AI and AI-assisted technologies in the writing process}
During the preparation of this work, the authors used OpenAI ChatGPT to
improve the language of the manuscript and to check the
consistency of notation. OpenAI ChatGPT and OpenAI Codex were used as coding assistants during the development of the simulation and optimization code. The tool was not used to determine the scientific conclusions of this work.

\section*{Data availability}
The Python implementation of the coupled synthesis--separation model, including the calibration of the chromatographic model and the scripts reproducing all figures of this manuscript, is openly available at \url{https://doi.org/10.5281/zenodo.XXXXXXX}. No experimental data were generated in this study; the chromatographic model is calibrated against the published data of \cite{Hartmann2025Chromatography}.

\bibliographystyle{cas-model2-names}
\bibliography{literature}

\end{document}

%% file: acronyms.tex
\newcommand*{\NewAcronym}[4][]{%
  \DeclareAcronym{#2}{%
    short = #3,
    long = #4,
    foreign-plural={}, 
    #1
  }
}
\NewAcronym{emom}{eMoM}{exact method of moments}
\NewAcronym{fpp}{FPP}{fixed point equation}
\NewAcronym{edm}{EDM}{equilibrium-dispersive model}
\NewAcronym{ntp}{NTP}{number of theoretical plates}
\NewAcronym{ode}{ODE}{ordinary differential equations}
\NewAcronym{pde}{PDE}{partial differential equations}
\NewAcronym{pbe}{PBE}{population balance equations}

\DeclareAcronym{psd}{
  short = PSD,
  long  = particle size distribution,
  short-plural = s,
  long-plural  = s
}
\NewAcronym{rsd}{RSD}{relative standard deviation}
\NewAcronym{sd}{SD}{standard deviation}
\NewAcronym{rtd}{RTD}{residence time distribution}
\NewAcronym{sec}{SEC}{size-exclusion chromatography}

%% file: definitions.tex
\newcommand*{\cC}{\mathcal{C}}

\newcommand*{\cG}{\mathcal{G}}

\newcommand*{\cN}{\mathcal{N}} 

\newcommand*{\cP}{\mathcal{P}}

\newcommand*{\cR}{\mathcal{R}}

\newcommand*{\cV}{\mathcal{V}}

\newcommand*{\cY}{\mathcal{Y}}

\newcommand*{\D}{\mathrm{d}}

\newcommand*{\R}{\mathbb{R}}

\newcommand{\growth}{\textit{growth}\xspace}
\newcommand{\nucleation}{\textit{nucleation}\xspace}

\newcommand{\cpre}{c_{\mathrm{pre}}} 
\newcommand{\ctot}{c_{\mathrm{mon}}} 

\newcommand{\pfun}{p_{\mathrm{in}}} 
\newcommand{\p}{p} 
\newcommand{\kr}{k_\mathrm{r}} 
\newcommand{\TC}{T_{\mathrm{c}}} 
\newcommand{\cGs}{\cG_{1}} 
\newcommand{\cGc}{\cG_{0}} 
\newcommand{\xn}{x_n} 
\newcommand{\ym}[1]{y_{#1}} 

\newcommand{\Lref}{L_{\mathrm{c,ref}}} 

\newcommand{\LS}{L_{\mathrm{s}}} 
\newcommand{\LC}{L_{\mathrm{c}}}

\newcommand{\muCref}{\mu_{\mathrm{c, ref}}} 

\newcommand{\muCana}{\mu_\mathrm{c,ana}}
\newcommand{\muCtr}{\mu_{\mathrm{c,tr}}}

\newcommand{\ueff}{u_\mathrm{eff}}
\newcommand{\Deff}{D_\mathrm{eff}}

\newcommand{\nfacone}{n_{\mathrm{fac},1}}
\newcommand{\nfactwo}{n_{\mathrm{fac},2}}
\newcommand{\gfac}{g_{\mathrm{fac}}}

\newcommand{\ntpref}{\mathrm{NTP}_\mathrm{ref}}
\newcommand{\alphaNTP}{\alpha_{\mathrm{NTP}}}

\newcommand{\qfrac}{q_{\mathrm{frac}}} 
\newcommand{\qCin} {q_{\mathrm{c,in}}} 
\newcommand{\qCout}{q_{\mathrm{c,out}}} 
\newcommand{\qC} {q_{\mathrm{c}}} 
\newcommand{\qCinit} {q_{\mathrm{c},0}}
\newcommand{\qS} {q_{\mathrm{s}}} 
\newcommand{\qSout}{q_{\mathrm{s,out}}} 

\newcommand{\VolS}{V_\mathrm{s}}

\newcommand{\rtdC}{f_\mathrm{c}} 
\newcommand{\rtdCcum}{F_\mathrm{c}} 
\newcommand{\rtdS}{f_\mathrm{s}} 

\newcommand{\sigmaCana}{\sigma_\mathrm{c,ana}}

\newcommand{\sigmaCref}{\sigma_{\mathrm{c, ref}}}
\newcommand{\HCcol}{H_{c,\mathrm{col}}}

\newcommand{\tend}{t_1}
\newcommand{\tstart}{t_0} 
\newcommand{\tinj}{t_{\mathrm{inj}}} 
\newcommand{\control}{\mathbf{v}}

\newcommand{\Vflow}{Q_\mathrm{s}} 
\newcommand{\xd}{x_\mathrm{\star}} 
\newcommand{\msd}{d} 
\newcommand{\rmsd}{d_{\mathrm{r}}} 

\newcommand{\ptop}{\mathcal{Q}_{\mathrm{frac}}}
\newcommand{\yield}{\cY}
\newcommand{\profit}{\cP}
\newcommand{\revenue}{\cR}
\newcommand{\cost}{\cC}
\newcommand{\costfac}{\gamma_\ell}

\newcommand{\dc}{\dot{c}}
\newcommand{\dcpre}{\dot{c}_\mathrm{pre}}
\newcommand{\dctot}{\dot{c}_\mathrm{mon}}
\newcommand{\dym}[1]{\dot{y}_{#1}}

%% file: 1_Introduction.tex

\section{Introduction}
\label{sec:intro}

Nanoparticles, i.e., particles with at least one dimension in the range of 1 to 100 nm (\cite{EC2022Nanomaterial}), exhibit properties that differ substantially from those of their bulk counterparts, for example with respect to their optical (\cite{seifert2024continuous, traor2024targeted, guo2023intrinsic}), electronic (\cite{parihar2018brief, snchezromate2021electrical}), and catalytic (\cite{lim2024nanoparticle, hoeven2021unlocking}) behavior.  These macroscopic properties are governed by the disperse properties of the particle ensemble, i.e., its size, shape, and composition distributions, towards which many of the relevant functionalities show a pronounced sensitivity, such that even comparatively small variations in, e.g., particle size can translate into appreciable shifts in the resulting property (\cite{salaheldin2017automated}). Consequently, narrow disperse property distributions are required in order to obtain reproducible and application-relevant functionality and thus the width of this distribution can be regarded as a measure of product quality.  \newline

Two complementary strategies may be used to establish such narrow distributions. The first is to control the synthesis itself such that particles with the targeted disperse properties are formed directly, which requires identifying suitable process conditions.  Since the relationships between process parameters and the resulting disperse properties are generally nonlinear and rarely accessible by intuition, this identification can be performed efficiently in a model-guided manner, using a process function that maps the process parameters onto the disperse properties of the resulting nanoparticles (\cite{traor2026inverse, wang2023population}). Approaches of this kind have been applied successfully to the design of a variety of nanoparticulate products, including the targeted adjustment of the size, composition, and optical response of mono and bimetallic systems (\cite{bleyer2025predictive, traor2024mechanistic, traor2024targeted, seifert2026bridging, forestiere2015inverse, christiansen2020inverse, kashiwagi2022highthroughput}). Where syntheses of nanoparticulate materials are operated continuously, which facilitates reproducible and scalable production, continuous flow reactors provide well-defined and adjustable reaction conditions and have been employed for the controlled synthesis of a broad range of nanoparticle systems (\cite{Biegel2023, seifert2025silver, vlkl2023gramscale, yoko2025precision, danjo2025roomtemperature, watanabe2025synthesis}). \newline

The second strategy is to apply a downstream classification step in which the distribution obtained from the synthesis is actively narrowed. On the laboratory scale, several such techniques are established (\cite{kowalczyk2011nanoseparations}), including centrifugation, which separates particles according to their sedimentation behavior (\cite{sun2009separation}), selective precipitation (\cite{segets2013quantitative, rezvani2025separation, guo2011shape}), field-flow fractionation (\cite{kammer2011separation, dubascoux2008optimisation}), dielectrophoretic separation (\cite{pesch2020review}), or liquid chromatography (\cite{gromotka2022classification}). Chromatographic classification in particular has been demonstrated for a broad spectrum of nanoparticulate systems, including spherical (\cite{gromotka2022classification, Hartmann2025Chromatography}) and anisotropic gold nanoparticles (\cite{hartmann2025characterization}), metal nanoclusters (\cite{gromotka2023green, supper2023isolation}), and carbon nanomaterials (\cite{s2019quantitative, hinterberger2019purification}), using size-exclusion and interaction based methods.
In addition, for chromatographic separations, process models are likewise available that predict the fractionation outcome as a function of the operating conditions and thereby permit a model-guided optimization analogous to that of the synthesis (\cite{SchmidtTraub2012, Snyder2009}). Both strategies can finally be combined by optimizing the whole process chain where a synthesis process is coupled together with a subsequent classification step that further reduces the width of the \ac{psd}. This third route can be expected to yield the narrowest distributions. \newline

For real-world applications, such a combined process chain must be realized not only on the laboratory scale but also at production-relevant throughput, which considerably restricts the set of admissible techniques, since batch processes and several of the aforementioned classification routes do not easily transfer to larger scales (\cite{panariello2020highly, kohsakowski2019effective, GROMOTKA2024_Coupling}). One combination that offers scalability is the coupling of a continuous flow nanoparticle synthesis with subsequent chromatographic classification. We have previously realized such a setup experimentally, where gold nanoparticles were produced in a continuous flow T-mixer setup, collected in an intermediate storage tank, and subsequently injected into a chromatographic column, with the throughputs of synthesis and classification adapted to one another (\cite{GROMOTKA2024_Coupling}). In that setup, \ac{sec} was used for the classification, during which nanoparticles are separated according to their hydrodynamic diameter (\cite{gromotka2022classification}). \newline

Although process models exist for both unit operations of such a process chain, the model-guided optimization of nanoparticulate systems has to date remained largely confined to single unit operations (\cite{nirschl2021autonomous, uhlemann2020modeling, traor2026inverse}). This is in contrast to general chemical engineering, where the flowsheet simulation of interconnected unit operations is well established for process design (\cite{sinnott2019chemical}) and end-to-end process models have been reported, for example, for continuous pharmaceutical manufacturing (\cite{nagy2021dynamic}). To the best of our knowledge, a coupled process chain model that links a nanoparticle synthesis to a downstream classification and which treats the dispersity of the particle ensemble as the quantity to be optimized across the entire chain has not yet been reported. Such a description within a single model is, however, the natural way to simulate and optimize this type of process chain, as it allows the particle dispersity to be shaped along the full chain rather than through the synthesis step alone, and thus permits trade-offs between the two unit operations to be resolved jointly rather than sequentially. \newline

\paragraph{Our contribution}
In this work, we present a coupled process model for continuous nanoparticle synthesis followed by classification using \ac{sec}. First, we connect the two unit operations of synthesis and separation by directly transferring the particle number distribution obtained at the synthesis reactor outlet to the inlet of the chromatographic model, thereby providing a consistent description of the complete process chain. Second, the \acl{emom}-based synthesis model is combined with analytical expressions for chromatographic transport, enabling efficient evaluations of the coupled process. Third, we formulate an economically driven optimization problem in which synthesis and separation controls are optimized simultaneously with respect to product quality, yield, and chromatographic cost. This formulation reveals process wide trade-offs that cannot generally be exploited when the two unit operations are designed sequentially.

Since the model is formulated for a prototypical gold-like model system rather than for a specific experimental material, the methodology is transferable to other dispersed particle systems and compatible downstream classification processes. It thereby provides a basis for extending model guided inverse design from individual unit operations to full nanoparticle process chains.

\paragraph{Outline of this work} In Section \ref{sec:model}, we start by introducing a mathematical model for particle synthesis governed by population balance equations and obtain computational results for a prototypical, gold-like system. Subsequently, we derive a mathematical model for particle separation in the equilibrium-dispersive model and again produce computational results for the gold-like system. Finally, the integrated process chain of coupled synthesis and separation is outlined. Section~\ref{sec:optimization} discusses the objective function for the optimized synthesis-separation process chain, taking into account the conflicting goals of obtaining a low-cost process which at the same time generates a high-quality product. Optimized results are presented and discussed for the gold-like model system. In Section~\ref{sec:conclusion} we conclude with a brief summary and potential directions for future research.  

%% file: 2_Simulation.tex


\section{Coupled synthesis and separation model}
\label{sec:model}
\subsection{Process setup}
\label{sec:processsetup}

In this work, the aim is to model an entire process chain, from the continuous synthesis of nanoparticles to the classification of the formed particles by \ac{sec}. Both unit operations are designed to act together rather than in sequence, such that the chain as a whole yields a particulate product of targeted mean diameter and dispersity. Throughout this work, the process chain is built around a generic gold-like nanoparticle system. While the synthesis step is not matched to a specific experimental counterpart, the assumed chemistry is that of gold, and the chromatographic model is calibrated against \ac{sec} experiments on gold nanoparticles (see \ref{sec:separation}). \\

\begin{figure}[pos=htbp]
    \centering
    \includegraphics[width=\textwidth]{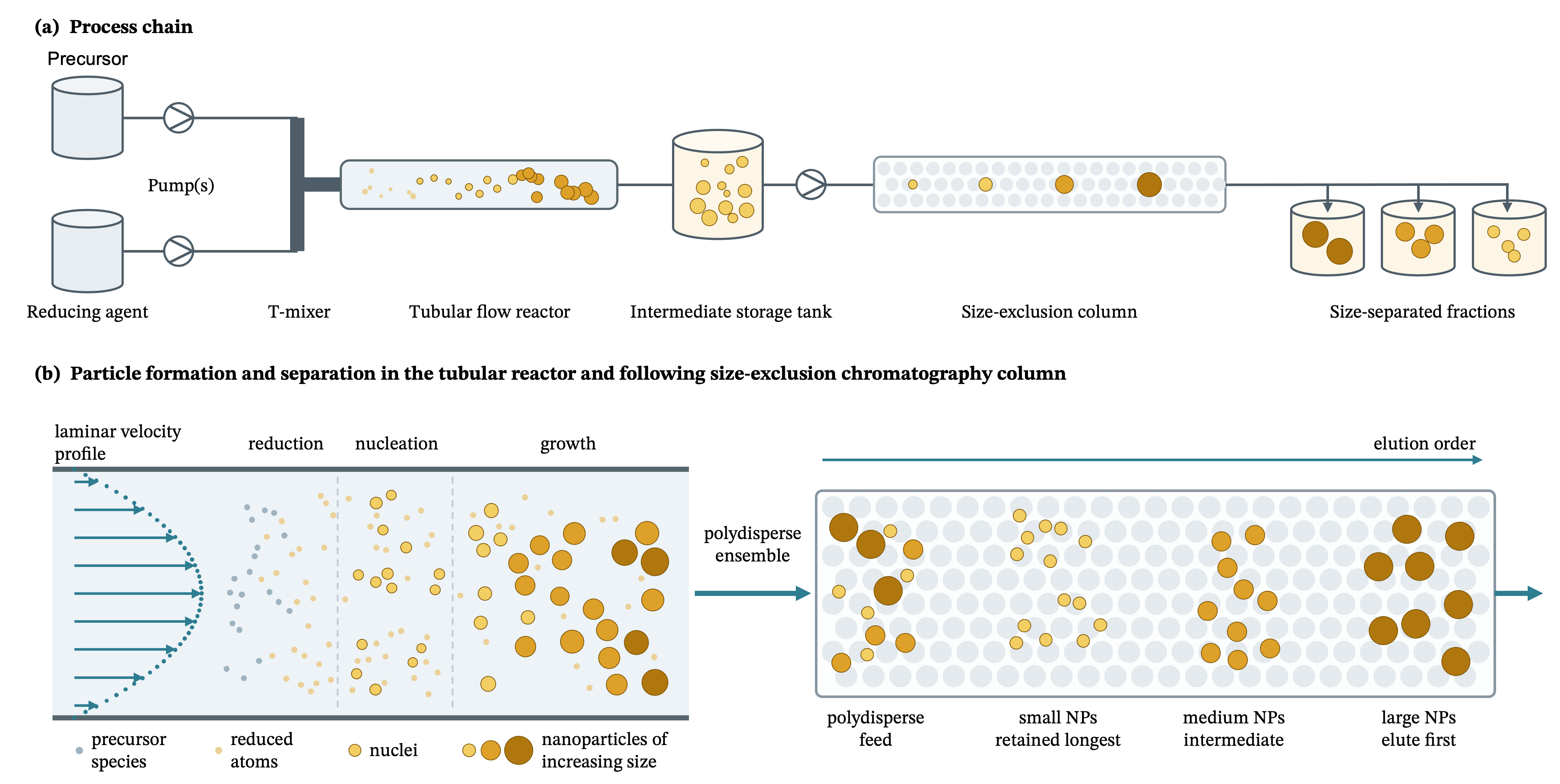}
    \caption{Process setup modeled in this work. a) Process chain starting with the mixing of precursor and reducing agent in a continuous T-mixer setup connected to a tubular reactor in which particle formation takes place. The product dispersion is collected in an intermediate storage tank, where the reaction is assumed to stop so that the present \ac{psd} is frozen. From there, the nanoparticles are periodically injected into the chromatographic column, which separates the particles according to their size. b) Axial section of the tubular reactor (left) and of the chromatographic column (right). The reactor shows the parabolic flow profile characteristic of the laminar flow regime. Particle formation proceeds along the axial coordinate and comprises the reduction of the precursor to zero-valent atoms, their nucleation and the subsequent growth to the final particle size. The resulting ensemble is disperse in size, which motivates its classification in the downstream chromatographic column, where the particles are separated along the packing according to their accessibility to the pore volume and therefore according to their size. Large nanoparticles thus exit the column first, whereas smaller nanoparticles are retained longer. }
    \label{fig:processsetup}
\end{figure}

\cref{fig:processsetup} shows the simulated process setup. In a T-mixer, a precursor solution is combined with a solution of the reducing agent. At sufficiently high flow rates, the energy coupled into the system generates vortices in the mixer and thus efficient mixing of the two streams, which can be characterized by the dimensionless Reynolds number (\cite{Biegel2023, schikarski2019influence, schikarski2019inflow}). Throughout this work we consider a regime in which mixing is complete on a timescale short compared with that of nucleation.  \newline
Particle formation takes place in the subsequent tubular reactor under isothermal conditions at room temperature. Reduction of the precursor releases monomers, which build up supersaturation, which in turn drives the formation of nuclei (nucleation) and their growth to the final particle size by monomer attachment to the surface of the existing nuclei. We disregard stabilizing ligands and assume the particles to be colloidally stable without surface active species. This simplification is introduced only for the present proof of principle. The effect of a stabilizer can in principle be included in the same framework. The particles leaving the reactor are collected in an intermediate storage tank, which is assumed to be perfectly mixed. We further assume that changes in particle morphology cease immediately upon entry into the tank, so that the \ac{psd} is effectively frozen. Experimentally, such quenching can be realized by a drop in temperature (\cite{wang2023population}), by the addition of thiolated ligands that render the particle surface inaccessible to further monomer attachment (\cite{Biegel2023, seifert2026bridging}), or by separating the reducing agent or the residual precursor, for example by means of membranes (\cite{sweeney2006rapid}). \newline
The particles leaving the reactor show a distribution in sizes for two reasons. First, nucleation and growth rates overlap in time, so that particles nucleated at different instants have grown for different durations. Second, the laminar flow profile in the tubular reactor imposes a \ac{rtd}. Under the assumption of complete segregation, volume elements travel through the reactor at their local radial velocity without exchanging material, so that elements close to the reactor wall remain in the reactor considerably longer than those close to the centerline. Both effects translate into a width of the final \ac{psd} as long as monomer is still available. Once the precursor is exhausted, additional residence time no longer changes the particle size, and the size differences accumulated within the reactor cease to develop further. \newline
In this process step, the \ac{psd} can therefore be influenced by the reaction time in the flow reactor, which at constant volumetric flow rate is set by the length of the tubular reactor, and by the reductive capacity of the system, which is set by the type of reducing agent and by its ratio to the precursor. For a given set of process conditions, a longer reactor yields larger particles and, provided that it is long enough for the precursor to be consumed, a narrower distribution, since all volume elements then approach the same final size distribution irrespective of their individual residence time. A higher reductive capacity accelerates the build-up of supersaturation and shortens the nucleation burst, so that a larger number of particles nucleate within a narrower time window, leading to a more uniform final \ac{psd}.
From the intermediate storage tank, the nanoparticles are fed into the chromatography column by periodic injections. Size exclusion chromatography is employed to further reduce the width of the \ac{psd}. Within the porous stationary phase, particles that show no enthalpic interaction with the stationary phase are separated according to their hydrodynamic diameter, because the fraction of the pore volume that is accessible to a particle decreases with increasing particle size. Since neither a ligand shell nor a hydration layer is considered here, we take the hydrodynamic diameter to be equal to the particle core diameter. Larger particles are thus excluded from a larger part of the pore volume, spend less time in the stagnant pore fluid and elute at a smaller retention volume, whereas smaller particles access a larger fraction of the pore volume and are retained longer. This process is approximated via the linear adsorption model (see \cref{sec:separation}). The recorded chromatogram therefore represents the \ac{psd} in inverted order, distorted by the nonlinear relation between particle size and retention volume and by band broadening along the column. A fictitious \ac{sec} detector is assumed to reproduce the eluting distribution without size dependent weighting, different to e.g., extinction based detectors. \newline
In this process step, the \ac{psd} can be influenced by the length of the \ac{sec} column and by fractionation. The resolution attainable in the column increases with column length (see \cref{sec:separation}), so that a longer column separates a given pair of particle sizes further. Building on this, the \ac{psd} can be narrowed further by fractionation, that is, by collecting only those particles that elute between a lower and an upper retention time and discarding the rest. Column length and the two fractionation bounds are thus the levers by which the size dispersity of the product may be controlled in the classification step. \\

With process models available for both the continuous synthesis and the size exclusion chromatography, including the parameters that govern the width of the \ac{psd}, a parameter set for the entire chain can be identified that yields particles of targeted diameter and dispersity. Since the levers of the two unit operations act on the same product property, and since the distribution leaving the reactor determines what the column has to work with, such a parameter set need not coincide with the combination of the individually optimal settings of the two unit operations. In the following, the two models are introduced in turn, beginning with the model of the continuous synthesis.
 
To clearly distinguish both unit operations, all synthesis-related quantities are denoted by the subscript s, whereas chromatography-related quantities are denoted by the subscript c. This convention is adopted throughout the remainder of the paper.

\subsection{Particle synthesis model}
\label{sec:synthesis}
\subsubsection{Model derivation}
As described above, we consider nanoparticle formation by nucleation and growth in a tubular reactor, taking into account the \ac{rtd} of the reactor, i.e., the distribution of the times over which individual fluid elements remain inside the reactor and the particles they carry continue to grow.

Following \citet{Pflug2020}, the synthesis model for the described scenario can be formulated as follows. The \ac{pbe} describe the evolution of a particle population governed by \nucleation and \growth,  both of which are driven by the supersaturation of the monomers. Its solution is the ac{psd} $\qS$ of spherical particles as a function of particle diameter $x$ and process time $t$. 
The precursor concentration $\cpre$ and the concentration of the monomers $\ctot$ are determined by the precursor mass inflow rate $\pfun \geq 0$ and evolve according to the first order reaction given by an \ac{ode}: 
\begin{equation}
\label{eq:firstorderreaction}
\begin{aligned}
\dcpre(t) &= \tfrac{1}{\VolS}\pfun(t) - \kr\, \cpre(t), 
& \cpre(0) &= 0, \\
\dctot(t) &= \kr\, \cpre(t), 
& \ctot(0) &= 0,
\end{aligned}
\end{equation}
where $\kr > 0$ is the reduction rate constant of the precursor species and $\VolS$ the considered volume of the volume element.

The growth kinetics are given by the concentration-dependent growth rate $\cGc$ and the size-dependent growth rate $\cGs$. The \ac{pbe} then reads:
\begin{equation}
\label{eq:populationbalanceequations}
\begin{aligned}
\partial_t \qS(t,x) + \partial_x \big( \cGc(c(t)) \cGs(x)\, \qS(t,x) \big) &= 0, \\
\qS(0,x) &= 0, \\
\cGc(c(t)) \cGs(\xn)\, \qS(t,\xn) &= \cN(c(t)),
\end{aligned}
\end{equation}
where $\cN$  describes the concentration-dependent nucleation rate and $\xn \in \R_{\geq 0}$ the nucleation size.  Lastly, $\qS(t,x)$ describes the ac{psd} over time and is defined as a number density distribution, wherein the particle diameter $x \in \R_{\geq\xn}$ and $t \geq 0$. The liquid-phase concentrations $c,\ctot$ and the ac{psd} are coupled through conservation of mass. Since $\qS$ denotes the ac{psd} with respect to particle size, its third moment determines the total
particle volume. Assuming spherical particles with density $\rho$, the
corresponding total particle mass at time $t$ is given by \(\tfrac{\rho \pi}{6} \int_{\xn}^{\infty} x^3 \qS(t,x)\, \D x\). Furthermore, for a given volume element $\VolS$ the free monomer concentration is given by 
\begin{equation}
\label{eq:c_analytical}
c(t) =  \ctot(t) - \tfrac{\rho \pi}{\VolS 6} \int_{\xn}^{\infty} x^3 \qS(t,x)\, \D x.
\end{equation}
This coupled \ac{pbe}--\ac{ode} system describes the evolution of the ac{psd} and the reactant concentrations within a single fluid element as a function of its residence time.

In a continuous flow reactor, however, different fluid elements experience
different residence times. Consequently, the observed outlet \ac{psd} is
obtained as an \ac{rtd}-weighted ac{psd}. We model
the reactor as a straight circular tube with ideal, steady, fully developed
Poiseuille flow and assume complete segregation, i.e., no radial or axial
exchange of precursor, monomer, or particles between fluid elements. For a
tube of length $\LS$ and inner diameter $D_{\mathrm{s}}$, operated
at the volumetric flow rate $\Vflow$, the reactor volume $V_R$ and mean residence time $\tau_s$ 
are
\begin{equation}
V_{\mathrm{R}}=\frac{\pi D_{\mathrm{s}}^2\LS}{4},
\qquad
\tau_{\mathrm{s}}=\frac{V_{\mathrm{R}}}{\Vflow}.
\label{eq:synthesis_mean_residence_time}
\end{equation}
The parabolic velocity profile has twice the cross-sectional mean velocity at
the centerline. Hence, the earliest possible exit time is
\begin{equation}
t_{\min}=\frac{\tau_{\mathrm{s}}}{2}
=\frac{\pi D_{\mathrm{s}}^2\LS}{8\Vflow}.
\label{eq:synthesis_min_residence_time}
\end{equation}
Flow-weighting the annular fluid elements gives the ideal laminar-tube exit-age
density \citep{Danckwerts1953,HartGuymerSonnenwaldStovin2016}
\begin{equation}
\label{eq:rtdS}
\rtdS(t)=
\begin{cases}
0, & 0\leq t<t_{\min},\\[1mm]
\displaystyle\frac{2t_{\min}^2}{t^3}, & t\geq t_{\min}.
\end{cases}
\end{equation}
Here, $\rtdS$ is a probability density with units of inverse time, rather than a
particle concentration or particle flow rate. It is nonnegative and normalized,
$\int_0^\infty\rtdS(t)\,\D t=1$, and its first moment, i.e., its mean residence time, is
$\int_0^\infty t\rtdS(t)\,\D t=2t_{\min}=\tau_{\mathrm{s}}$.

The reactor outlet \ac{psd} is therefore given by
\begin{align}
\qSout(x)=\int_0^\infty \qS(\tau,x)\,\rtdS(\tau)\,\D\tau.
\label{eq:qSout}
\end{align}

The formulation above defines the synthesis model at the level of an individual fluid element and its \ac{rtd}-weighted reactor outlet. For its efficient numerical evaluation, we next derive a reduced representation of the coupled PBE--\ac{ode} system.

\paragraph{Solution approach via the \ac{emom}}
Solving the precursor and total-species balances in
\eqref{eq:firstorderreaction} and substituting their analytical solutions
into the mass balance \eqref{eq:c_analytical} yields
\begin{equation}
\label{eq:c}
c(t)
=\frac{1}{\VolS} \Big(
\int_0^t
\bigl(1-e^{-\kr(t-\tau)}\bigr)\pfun(\tau)\,\D\tau
-
\tfrac{\rho\pi}{6}
\int_{\xn}^{\infty}
x^3\qS(t,x)\,\D x \Big),
\end{equation}
where $\VolS >0$ denotes the volume of a volume-element and $\pfun$ the precursor mass inflow. Since $\qS$ depends on the concentration through the nucleation and
growth kinetics, \eqref{eq:c} defines an implicit relation for $c$. To evaluate the coupled population-balance and mass-balance system
efficiently, we employ the \ac{emom} developed by \citet{Pflug2020}.

For the synthesis process considered here, we specialize the general
framework to reaction-limited growth, i.e., $\cGs\equiv1$. We further assume an impulsive precursor input at the inlet,
\begin{equation}
\label{eq:pin}
    \pfun(t)=\p\,\delta(t),
\end{equation}
representing the instantaneous supply of the entire precursor amount to each volume element at $t=0$. Under these
assumptions, the \ac{pbe} can be represented by a closed system of ordinary
differential equations for the concentration and selected moments (see e.g. \cite{Bnsch2023}), while
the full \ac{psd} is subsequently reconstructed from the \ac{emom}
solution. This closed system can be solved efficiently in practice. Following the derivation of \citet{Pflug2020}, and using the \ac{pbe}, integration by parts, and the boundary condition at the nucleation size, the derivative of the $m$-th moment is obtained as
\begin{align*}
\dym{m}(t)
&= \frac{d}{dt} \int_{\xn}^{\infty} x^m \qS(t,x)\,\D x \\
&= - \int_{\xn}^{\infty} x^m \partial_x\bigl( \cGc(c(t))\, \qS(t,x) \bigr)\,\D x\\
&= \xn^m \cN(c(t)) + m \cGc(c(t)) \ym{m-1}(t).
\end{align*}
Using the mass balance together with the impulsive precursor inflow \eqref{eq:pin}
the time derivative of the free monomer concentration is obtained as
\begin{equation*}
\dc(t)
=\frac{1}{\VolS}\Big(
\kr \p e^{-\kr t}
-
\frac{\rho\pi}{6}
\dym{3}(t)\Big).
\end{equation*}
\paragraph{Reformulation} In summary, the overall model for the continuous flow reactor synthesis is given by

\begin{subequations}\label{eq:synthesismodel}
\begin{align}
\dc(t)
&=
\tfrac{1}{\VolS}
\left(
\kr \p e^{-\kr t}
-
\tfrac{\rho\pi}{6}
\bigl(
\xn^3 \cN(c(t))
+
3 \cGc(c(t)) \ym{2}(t)
\bigr)
\right),
\label{eq:ODEc}\\
\dym{2}(t)
&=
\xn^2 \cN(c(t))
+ 2 \cGc(c(t)) \ym{1}(t),
\label{eq:ODEy2}\\
\dym{1}(t)
&=
\xn \cN(c(t))
+ \cGc(c(t)) \ym{0}(t),
\label{eq:ODEy1}\\
\dym{0}(t)
&=
\cN(c(t)),
\label{eq:ODEy0}
\\[1ex]
\qS(t,x(\tau,t))
&= \cN(c(\tau))(\cGc(c(\tau)))^{-1},
\label{eq:q}
\\
x(\tau,t)
&= \xn + \int_{\tau}^{t} \cGc(c(s))\, \D s \label{eq:x}
\end{align}
\end{subequations}
for each $t \geq 0$ and with initial values \(c(0)=\ym{0}(0)=\ym{1}(0)=\ym{2}(0)=0\).

Here, \eqref{eq:ODEc}-\eqref{eq:ODEy0} define a closed system of ordinary differential equations for the free monomer concentration $c$ and the moments $\ym{0}$, $\ym{1}$, and $\ym{2}$. Furthermore, \eqref{eq:q} and \eqref{eq:x} provide the formulas for the \ac{psd} and particle size $x(\tau,t)$ at time $t$, obtained via the \ac{emom} under the assumption $\cGs \equiv 1$, see \cite[Thm.~2.2]{Pflug2020} based on \cite[Thm.~2.20]{KEIMER20174023} and \cite[Thm.~3.1]{NonlocalBounded2018}.

Finally, \eqref{eq:qSout} models the outlet \ac{psd} from the continuous flow reactor using the analytical Poiseuille-flow \ac{rtd} given in \eqref{eq:rtdS}. The ODE system \eqref{eq:ODEc}-\eqref{eq:ODEy0} is solved numerically to obtain $c$, $\ym{0}, \ym{1}, \ym{2}$. Furthermore, the solution of the \ac{pbe} $\qS(x,t)$ as well as the outlet ac{psd} $\qSout$ is computed via a consistent numerical approximation scheme based on the same discretization. 
The process models were implemented in Python using NumPy and SciPy.
The synthesis moment equations were integrated with the implicit Radau method provided by \texttt{scipy.integrate.solve\_ivp}. The particle size  distribution was subsequently reconstructed from the characteristic  relations of the \ac{emom}, and the reactor outlet distribution was obtained by integration over the synthesis \ac{rtd}.

\subsubsection{Forward simulation of synthesis products}
\label{sec:forward_sim_syn}

With a model for nanoparticle formation in the continuous flow reactor at hand, we now study its dynamics, in particular the response of the process outcome to changes in the process parameters. A mechanistic model such as the one derived here should reproduce physically meaningful trends in the process outcome. Since the model is not calibrated against an experimental counterpart in this work, we first identified a reasonable set of process parameters assuming the chemistry of gold.

We consider the growth and nucleation rate functions according to \cite{Thanh2014,Mersmann2001}
\begin{equation}
\label{eq:growth_instance}
\cGc(c)
=
\begin{cases}
    \gfac(S(c)-1), & S(c)>1,\\

    0,      & \text{otherwise},
\end{cases}
\end{equation}
and
\begin{equation}
\label{eq:nucleation_instance}
\cN(c)
=
\begin{cases}
\nfacone\, S(c)^{7/3}
\exp\!\left(
-\dfrac{\nfactwo}{(\ln (S(c)))^2}
\right),
& S(c)>1, \\[1ex]
  0,      & \text{otherwise},
\end{cases}
\end{equation}
which depend on the supersaturation $S(c)=c/c_\mathrm{eq}$ of the free monomer concentration $c$, with the monomer solubility $c_\mathrm{eq}=10^{-12}\,\mathrm{mol\,L^{-1}}$ (see \cref{tab:synth_parameters}). Here, $\gfac$ denotes the characteristic growth-rate scale with units of $\mathrm{m\,s^{-1}}$, whereas $\nfacone$ denotes the nucleation rate prefactor with units of $\mathrm{particles\,s^{-1}}$ and $\nfactwo$ is a dimensionless parameter controlling the supersaturation dependence of the nucleation rate. All parameters required for the synthesis simulation and their chosen values are summarized in Table ~\ref{tab:synth_parameters}.

\begin{table}[htbp]
\centering
\caption{Reference parameters for the nanoparticle synthesis model.}
\label{tab:synth_parameters}
\begin{tabular}{lccc}
\hline
\textbf{Parameter} & \textbf{Symbol} & \textbf{Value} & \textbf{SI-Unit}\\
\hline
Nucleation size & $\xn$ & $1.000\times10^{-9}$ & $\mathrm{m}$\\
Nucleation prefactor & $\nfacone$ & $4.500\times10^{-4}$ & $\mathrm{s^{-1}}$\\
Nucleation factor & $\nfactwo$ & $6.000\times10^{3}$ & -\\
Growth prefactor & $\gfac$ & $1.000\times10^{-17}$& $\mathrm{m\,s^{-1}}$\\
Particle density & $\rho$ & $1.932\times10^{7}$ & $\mathrm{g\,m^{-3}}$\\
Precursor mass & $\p$ & $2.500\times10^{-11}$ & $\mathrm{g}$\\
Volume element & $\VolS$ & $1.000 \times 10^{-12}$& $\mathrm{m^3}$\\
Equilibrium concentration & $c_{\mathrm{eq}}$ & $1.970\times10^{-7}$ & $\mathrm{g\,m^{-3}}$\\
Reactor diameter & $D_{\mathrm{s}}$ & $2.000\times10^{-3}$ & $\mathrm{m}$\\
Volumetric flow rate & $\Vflow$ & $1.667\times10^{-8}$& $\mathrm{m^3\,s^{-1}}$\\
Reduction rate constant & $\kr$ & $3.333 \times 10^{-4}$ -- $3.333 \times 10^{-3}$& $\mathrm{s^{-1}}$\\
Reactor length & $\LS$ & $2.000$--$25.000$ & $\mathrm{m}$\\
\hline
\end{tabular}
\end{table}

In general, the reactor outlet \ac{psd} is obtained by weighting the time-dependent synthesis ac{psd} with the reactor \ac{rtd}, see \eqref{eq:qSout}. Holding the remaining process parameters constant, the \ac{psd} at the reactor outlet $\qSout$ can be tuned by varying either the reduction rate constant $\kr$ or the reactor length $\LS$. A change in $\kr$ physically corresponds to a stronger reducing agent or to a higher concentration of the existing one, whereas $\LS$ determines the stage of the reaction reached at the outlet. A shorter reactor delivers smaller nanoparticles early in their growth phase, together with a wider \ac{psd}, since the parabolic flow profile still has a pronounced effect, see \ref{sec:processsetup}.

In the computational study, we therefore investigate the individual effects of both variables on the resulting \ac{psd}. To this end, we vary the reduction rate constant $\kr$ at a fixed reactor length $\LS$, and subsequently vary $\LS$ while keeping $\kr$ fixed.

\begin{figure}[pos=htbp]
    \centering
    \includegraphics[width=\textwidth]{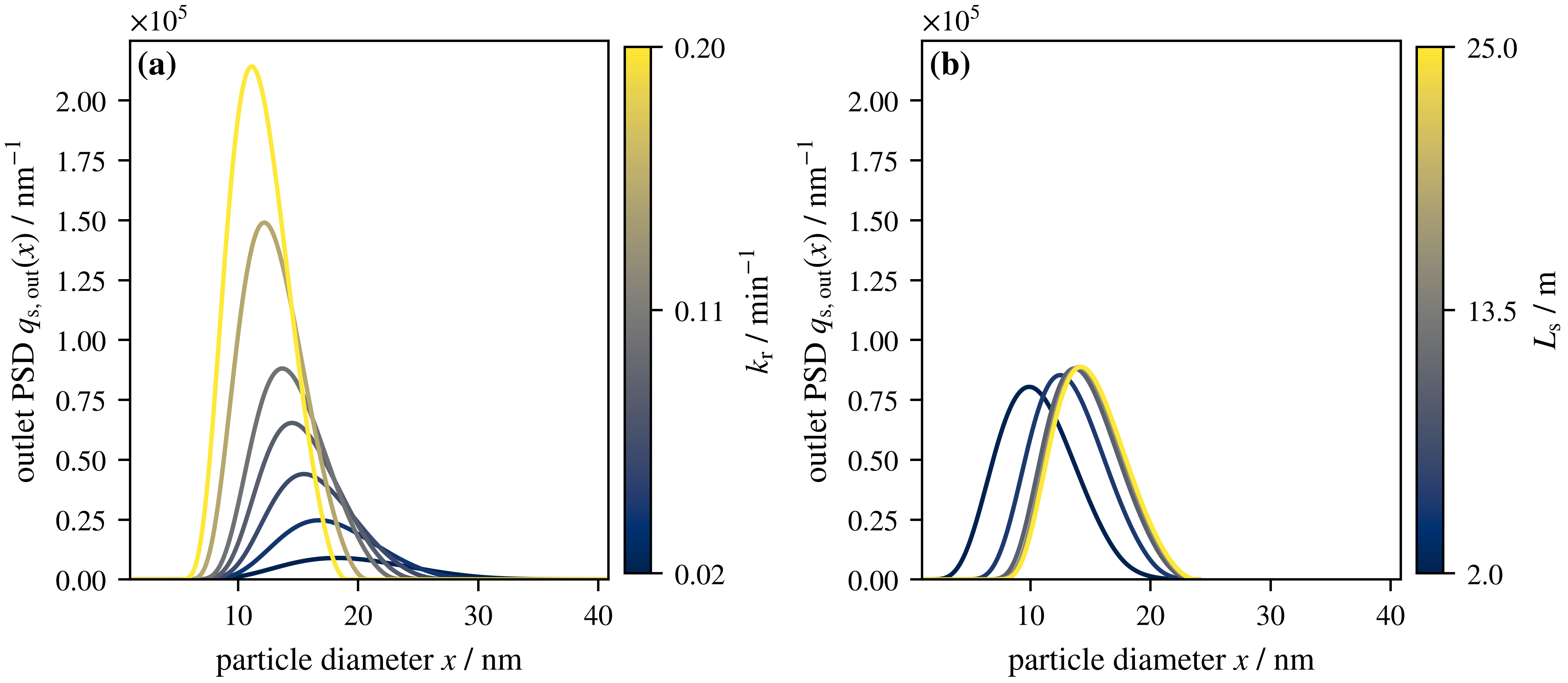}
    \caption{Influence of the synthesis controls on the reactor outlet \ac{psd}~$\qSout$: (a) variation of the reduction rate constant $\kr$ at $\LS=10\,\mathrm{m}$ and (b) variation of the synthesis reactor length $\LS$ at $\kr=0.1\,\mathrm{min}^{-1}$.}
    \label{fig:psd_k_L}
\end{figure}

\paragraph{Variation of the reduction rate constant} Figure~\ref{fig:psd_k_L}a shows the number weighted reactor outlet \acp{psd} obtained for different reduction reaction rate constants $\kr$ between $0.02$ and $0.20\,\mathrm{min}^{-1}$ at a constant reactor length $\LS=10\,\mathrm{m}$. The rate constant determines how rapidly precursor is converted into monomers (see \ref{eq:firstorderreaction}), and therefore governs the supersaturation profile over time. Since both nucleation and particle growth depend on the supersaturation (see \ref{eq:growth_instance} and \ref{eq:nucleation_instance}), varying $\kr$ changes their relative contributions over time and consequently affects the mean size and the width of the reactor outlet \ac{psd}. High values of $\kr$ shift the \ac{psd} towards smaller mean particle sizes and narrow its width as the rapid rise in supersaturation produces a large number of nuclei in a short burst, and the accompanying rapid consumption of precursor prevents the supersaturation from returning to levels that would support significant further nucleation at later times. The existing nuclei then only grow, which, under the assumption of size-independent integration-limited growth and ideal mixing, does not broaden the \ac{psd} further.

\paragraph{Variation of the reactor length} Figure~\ref{fig:psd_k_L}b shows the reactor outlet \acp{psd} for reactor lengths $\LS$ between $2$ and $25\,\mathrm{m}$ at a constant reduction rate constant of $\kr=0.1\,\mathrm{min}^{-1}$. The reactor length determines the synthesis \ac{rtd}. Increasing $\LS$ shifts the \ac{rtd} towards longer process times, allowing precursor reduction and particle growth to proceed further before exiting the reactor. Consequently, the outlet \ac{psd} shifts towards larger particle diameters and becomes slightly narrower. Unlike $\kr$, $\LS$ does not change the intrinsic kinetic laws, but determines how far the same synthesis process has progressed at the reactor outlet.
\newline
Overall, the figure demonstrates that the synthesis stage already provides key control parameters for the \ac{psd} and, thus, for product quality. At the same time, synthesis alone does not generally yield the desired combination of mean particle size, narrow dispersity, and high product yield. A short reactor may deliver nanoparticles of the targeted diameter, but at the expense of a comparably wide \ac{psd} and large amounts of unreacted precursor.  A high reduction rate constant, in turn, may yield narrowly distributed particles without allowing the mean particle size to be tuned to the desired value at the same time. This motivates the subsequent chromatographic classification of the synthesis product and, ultimately, the joint optimization of both process stages towards a targeted \ac{psd}.

\subsection{Particle separation in the equilibrium-dispersive model}
\label{sec:separation}
\subsubsection{Model derivation}

The starting point for the separation model is a commonly used \ac{edm}, see \citet{SchmidtTraub2012}. This model describes the evolution of the ac{psd} $\qC$ over the residence time $t \in [0,\TC]$ and the spatial coordinate $z \in [0,L]$ along the chromatographic column:

\begin{equation}
\begin{aligned}
\partial_t \qC(t,x,z) + A\,\partial_t a\!\left(\qC(t,x,z)\right) + u(x)\,\partial_z\qC(t,x,z)&=D(x)\,\partial_z^2\qC(t,x,z),\\
u(x)\qC(t,x,0)-D(x)\partial_z\qC(t,x,0)&=u(x)\qCinit(t,x),\\
\partial_z\qC(t,x,L)&=0\\
\qC(0,x,z)&=0.
\end{aligned}
\label{eq:equilibriumdispersivemodel}
\end{equation}

Here, $A$ denotes the phase ratio, $u$ the velocity, $D$ the dispersion coefficient, and $a$ the adsorption isotherm. The inflow is described by a Danckwerts boundary condition at $z=0$ (\cite{Danckwerts1953}), the outlet by a free outflow condition at $z=L$, and the column is assumed to be initially empty. In our work, we consider \ac{sec} which is often modeled via a linear function for the adsorption isotherm. The adsorption isotherm is then given by $a(\qC) = K(x) \qC$ with a size dependent distribution coefficient $K$. In this case, the adsorption term simplifies to $A\,\partial_t a(\qC) = A K(x)\, \partial_t \qC$.
Consequently, the \ac{edm} \eqref{eq:equilibriumdispersivemodel}
is linear in $\qC$. Moreover, since no derivatives with respect to the particle size
$x$ appear, every particle size can be considered independently. Introducing effective velocity and diffusion coefficients, i.e., 
\begin{align*}
\ueff := \tfrac{u(x)}{1+AK(x)}, \quad \text{and} \quad \Deff := \tfrac{D(x)}{1+AK(x)}     
\end{align*} 
\eqref{eq:equilibriumdispersivemodel} can be rewritten as
\begin{equation}
\label{eq:edm_linear}
    \partial_t \qC(t,x,z)
    +\ueff(x)\partial_z \qC(t,x,z)
    =
    \Deff(x)\partial_{z}^2\qC(t,x,z).
\end{equation}

To eliminate the advective transport, we introduce the moving coordinate system by substituting $z = \tilde{z} + \ueff(x)t $

\begin{equation*}
    P(t,x,\tilde{z})
    :=
    \qC\!\left(
        t,x,
        z+\ueff t
    \right).
\end{equation*}
Applying the chain rule gives
\(
    \partial_t P
    =
    \partial_t \qC
    +
    \ueff\partial_{z} \qC,
    \partial_{\tilde{z}}^2P
    =
    \partial_{z}^2\qC.
\)
Hence, \eqref{eq:edm_linear} transforms into
\(\partial_t P
    =
    \Deff(x)\partial_{\tilde{z}}^2P,\)
which has the form of the well-known heat equation that describes how heat distributes over a given region over time.
Thus, in the moving coordinate system, the chromatographic transport equation
reduces to the heat equation. In other words, the transport of nanoparticles through the chromatographic column can be mathematically described as heat or mass transport. The fundamental solution of the heat equation is the Gaussian distribution, see e.g., \citet{Evans2010}:
\begin{align}
\label{eq:heat_kernel_chrom}
P(t,x,\cdot)
\sim
\mathcal N\!\left(
0,\,
2\Deff(x) t
\right).
\end{align}

Transforming back to the original coordinates, the \ac{edm} for linear chromatography, as discussed in \cite{SchmidtTraub2012,Snyder2009}, can be approximated by a Gaussian:
\begin{align}
\label{eq:adv_disp_solution_chrom}
    \qC(t,x,z)
    =
    \tfrac{1}{
    \sqrt{4\pi \Deff(x)t}}
    \exp\left(
        -
        \tfrac{\left(z-\ueff(x)t\right)^2}{ 4\Deff(x)t}
    \right).
\end{align}
This describes the advection and dispersion of the ac{psd} $\qC$ along the column axis originating from a spatial (dirac) pulse at the column inlet $\qC(0,x,z)=\delta(z)$ while neglecting boundary condition effects.
We denote the cumulative fraction of particles that has passed the position $z=L$ by $\rtdCcum(t,x;L)$ and the corresponding residence time density by $\rtdC(t,x;L)=\partial_t\rtdCcum(t,x;L)$:
\begin{align}
    \rtdCcum(t,x;L) &= \int_{L}^{\infty} \qC(t,x,z)\,\D z = \Phi\left(\tfrac{t\ueff(x)-L}{\sqrt{2t\Deff(x)}}\right), \notag\\
    \rtdC(t,x;L) &= \partial_t\rtdCcum(t,x;L) = \tfrac{t\ueff(x)+L}{4\sqrt{\pi\Deff(x)t^3}}\exp\left(-\tfrac{\left(t\ueff(x)-L\right)^2}{4\Deff(x)t}\right) \label{eq:f_c_analytical}
\end{align}
Thus, for each particle size $x$, the chromatographic column is characterized by a size-dependent \ac{rtd} determined by the effective transport velocity, the axial dispersion coefficient, and the column length.

\paragraph{Calibration}
Although the \ac{edm} provides an analytical approximation for the outlet concentration profile and the corresponding \ac{rtd} $\rtdC$ (see \eqref{eq:f_c_analytical}), calibration of chromatographic processes is typically performed in terms of retention times measured at the column outlet rather than the transport parameters $u(x), K$ and, $D(x)$ as  the latter cannot easily be measured 
experimentally \citep{Snyder2009}. \\

The analytical mean retention time and its variance are given by:
\begin{align}
\label{eq:mu_derivation}
\muCana(x;L)&=\tfrac{L}{\ueff(x)}+\tfrac{\Deff(x)}{\ueff(x)^2},\\ \notag
\sigmaCana(x;L)^2&=\tfrac{2\Deff(x)L}{\ueff(x)^3}+\tfrac{5\Deff(x)^2}{\ueff(x)^4}.
\end{align}

Instead of determining the individual transport parameters explicitly, we calibrate the effective transport parameters directly from experimental data. Let $\muCref(x)$ denote the measured mean retention time for a reference column of length $\Lref$ and let $\ntpref(x)$ denote the experimentally determined number of theoretical plates \citep{Snyder2009}, given by
\begin{equation}
\ntpref(x)=\left(\tfrac{\muCref(x)}{\sigmaCref(x)}\right)^2.
\label{eq:NTPbeta}
\end{equation}

Moreover, defining
\[
\beta(x):=\tfrac{1-\ntpref(x)+\sqrt{\ntpref(x)\left(\ntpref(x)+3\right)}}{5\ntpref(x)-1},
\]
the calibrated effective transport parameters are given by
\begin{align*}
\Deff(x)&=\frac{\Lref^2\beta(x)\left(1+\beta(x)\right)}{\muCref(x)},\\
\ueff(x)&=\frac{\Lref\left(1+\beta(x)\right)}{\muCref(x)}.
\end{align*}

Substituting these expressions into \eqref{eq:f_c_analytical} yields the calibrated cumulative \ac{rtd} for $t>0$ based on a reference column:
\begin{equation*}
\rtdCcum(t,x;L)=\Phi\left(\frac{\frac{\Lref(1+\beta(x))}{\muCref(x)}t-L}{\sqrt{2\frac{\Lref^2\beta(x)(1+\beta(x))}{\muCref(x)}t}}\right), \qquad t>0.
\end{equation*}
Here, $\Phi$ denotes the cumulative distribution function of the standard normal distribution:
\begin{equation*}
\Phi(z)
=
\frac{1}{\sqrt{2\pi}}
\int_{-\infty}^{z}
\exp\left(-\frac{s^2}{2}\right)\,\mathrm ds \qquad z\in\mathbb{R}.
\end{equation*}

\paragraph{Fractionation}
The chromatographic model derived above provides the cumulative intrinsic \acl{rtd} $\rtdCcum(t,x;L)$ and its density $\rtdC(t,x;L)=\partial_t\rtdCcum(t,x;L)$ for particles of size $x$. For an impulsive injection, the time-resolved outlet distribution is obtained by multiplying $\rtdC$ with the inlet ac{psd} $\qCin(x)$. In practice, however, particles are injected over a finite time interval. We therefore assume a rectangular injection pulse of duration $\tinj$. Owing to the linearity of the \acl{edm}, the outlet \ac{psd} is obtained by convolution of the injection profile with the intrinsic $\rtdC$:
\begin{align}
\qCout(t,x)&=\frac{\qCin(x)}{\tinj}\int_0^{\tinj}\rtdC(t-\tau,x;L)\,\mathrm d\tau \notag\\
&=\frac{\qCin(x)}{\tinj}\left[\rtdCcum(t,x;L)-\rtdCcum(t-\tinj,x;L)\right].
\label{eq:qCout}
\end{align}

The number based \acl{psd} collected within the fractionation interval $[\tstart,\tend]$ is therefore given by
\begin{equation}
\label{eq:qfrac}
\qfrac(x)=\int_{\tstart}^{\tend}\qCout(t,x)\,\mathrm dt,
\end{equation}
which represents the amount of particles of size $x$ recovered in the selected fraction.

For efficient numerical evaluation, this integral can be expressed analytically. Introducing
\begin{equation}
\label{eq:mutr}
\muCtr(x;L)=\frac{L}{\ueff(x)}, \qquad N(x;L)=\frac{L\ueff(x)}{2\Deff(x)},
\end{equation}
an antiderivative of $\rtdCcum$ is given by
\begin{align*}
\HCcol(t,x;L)
={}&
\left(t-\muCtr(x;L)\right)
\Phi\!\left(
\left(t-\muCtr(x;L)\right)
\sqrt{\tfrac{N(x;L)}{\muCtr(x;L)t}}
\right)
\\
&\quad
+\sqrt{\tfrac{\muCtr(x;L)t}{N(x;L)}}\,
\phi\!\left(
\left(t-\muCtr(x;L)\right)
\sqrt{\tfrac{N(x;L)}{\muCtr(x;L)t}}
\right)
\\
&\quad
-\tfrac{\muCtr(x;L)}{2N(x;L)}
\Bigg(
\Phi\!\left(
\left(t-\muCtr(x;L)\right)
\sqrt{\tfrac{N(x;L)}{\muCtr(x;L)t}}
\right)
\\
&\qquad\qquad
-\exp\!\left(2N(x;L)\right)
\Phi\!\left(
-\left(t+\muCtr(x;L)\right)
\sqrt{\tfrac{N(x;L)}{\muCtr(x;L)t}}
\right)
\Bigg),
\end{align*}
for $t>0$, and $\HCcol(t,x;L)=0$ for $t\le0$. Consequently,
\begin{align}
\label{eq:qfracAnalytical}
\qfrac(x)=\tfrac{\qCin(x)}{\tinj}\Big(&\HCcol(\tend,x;L)-\HCcol(\tstart,x;L) \notag \\
&-\HCcol(\tend-\tinj,x;L)+\HCcol(\tstart-\tinj,x;L)\Big).
\end{align}
\subsubsection{Forward simulation of particle separation}
\label{sec:forward_sim_chrom}

Having established the model for particle separation, we now study its dynamics, analogously to the synthesis model in \ref{sec:forward_sim_syn}. A rectangular injection pulse of duration $\tinj$ is assumed. The size-dependent reference mean retention time is fitted to experimental \ac{sec} data for gold nanoparticles obtained with a $30\,\mathrm{cm}$ column with $100\,\mathrm{nm}$ nominal pore size at a flow rate of $0.5\,\mathrm{mL\,min^{-1}}$ by \citet{Hartmann2025Chromatography}. The measured retention volumes are converted to retention times and approximated by a third-order polynomial
\begin{equation}
\label{eq:muC_approx}
\muCref(x)
=
p_1x^3+p_2x^2+p_3x+p_4,
\end{equation}
with coefficients
\begin{equation*}
\label{eq:muC_coefficients}
\begin{aligned}
p_1&=-1.3728\times10^{24}\,\mathrm{s\,m^{-3}},&
p_2&=2.3765\times10^{17}\,\mathrm{s\,m^{-2}},\\
p_3&=-1.9971\times10^{10}\,\mathrm{s\,m^{-1}},&
p_4&=1.3922\times10^{3}\,\mathrm{s}.
\end{aligned}
\end{equation*}
The size-dependent reference plate number $\ntpref(x)$ is likewise approximated from the experimental data reported by \citet{Hartmann2025Chromatography}. For particle diameters below $8\,\mathrm{nm}$, the value at the lower boundary of the calibrated range is used. Above $50\,\mathrm{nm}$, the value at the upper boundary is retained and the evaluation is flagged as being outside the calibrated separation range. To represent different chromatographic separation capabilities, we scale the reference curve according to
\begin{equation}
\label{eq:ntp_scaling}
\ntpref^{\alpha}(x)=\alphaNTP\,\ntpref(x), \qquad \alphaNTP>0.
\end{equation}
The reference simulation uses $\alphaNTP=1$. In other scaling scenarios, $\ntpref^{\alpha}(x)$ replaces $\ntpref(x)$ in the definition of $\beta(x)$ \eqref{eq:NTPbeta}. The other required parameters for the simulation of the separation process are summarized in Table~\ref{tab:chrom_parameters}.
\begin{table}[htbp]
\centering
\renewcommand{\arraystretch}{1.2}
\caption{Reference parameters for the chromatographic separation model.}
\label{tab:chrom_parameters}
\begin{tabular}{lccc}
\hline
\textbf{Parameter} & \textbf{Symbol} & \textbf{Value} & \textbf{SI-Unit}\\
\hline
Reference column length       & $\Lref$ & $0.3$& m\\
Injection duration            & $\tinj$ & $1$& s\\
Reference plate-number curve & $\ntpref(x)$ & calibrated data fit  & --\\
Reference mean retention time & $\muCref(x)$
                                      & \eqref{eq:muC_approx} & s\\
\hline
\end{tabular}
\end{table}

\begin{figure}[pos=htbp]
    \centering
    \includegraphics[width=\textwidth]{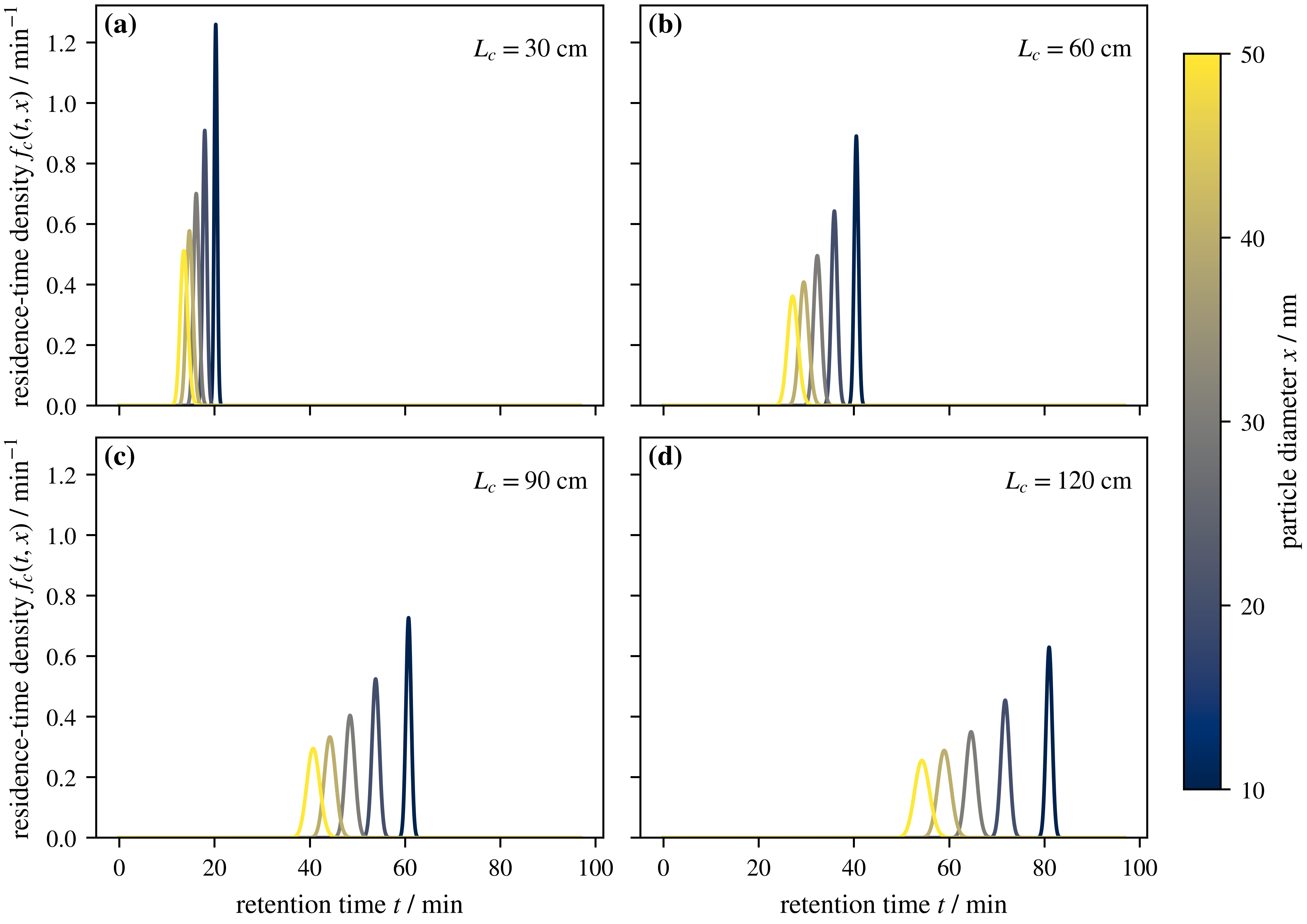}
    \caption{Chromatographic \acp{rtd} for gold nanoparticles with a diameter of 10-50 nm on size-exclusion columns of length 30-120 cm. The size dependence of the mean retention time and the column length enables chromatographic classification, while the peak width reflects band broadening, i.e., the temporal dispersion of particles of the same size at the column outlet.}
    \label{fig:fc}
\end{figure}

\begin{figure}[pos=htbp]
    \centering
    \includegraphics[width=\textwidth]{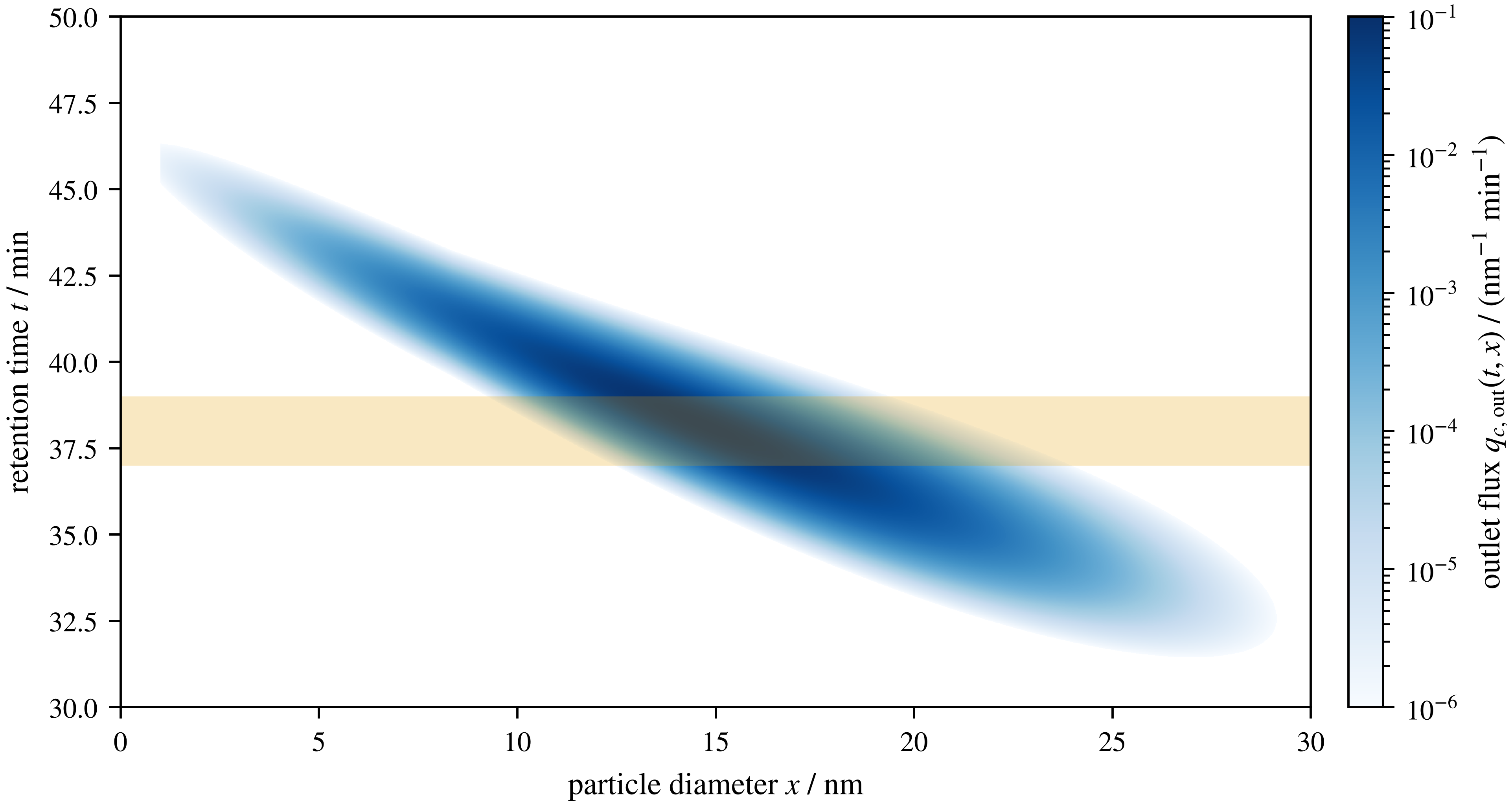}
    \caption{Time- and size-depended particle flux $\qCout$ at the outlet of the chromatographic column for a synthetic Gaussian inlet particle-number distribution. Color denotes the outlet flux in $(\mathrm{particles \;nm^{-1} \; min^{-1}})$ on a logarithmic scale. The orange band delimits the fractionation interval $([\tstart,\tend]=[37,39]\;\mathrm{min})$.}
    \label{fig:heatmap_qout}
\end{figure}

\Cref{fig:fc} shows the calculated chromatographic $\rtdC$ for monodisperse samples with particle sizes of 10-50\,nm assuming an impulse injection for column lengths of 30, 60, 90, and 120\,cm. For a given particle size, the \ac{rtd} is governed by two effects of the column length. First, the mean retention time given by \eqref{eq:mu_derivation} scales linearly with the column length, i.e., particles of the same size require twice the time to pass a column of twice the length. Second, the width of the \ac{rtd} scales with the square root of the column length, so that the \ac{rtd} broadens with increasing column length. In addition to these dominant scaling relations, both the mean retention time and its variance contain a column length independent correction arising from axial dispersion, see \eqref{eq:mu_derivation}. Superimposed on this is the size dependence, i.e., within the separation range of the given column, the mean retention time decreases monotonically with particle size according to the fitted retention
time relation, so that larger particles pass through the column faster than smaller ones (see \cref{sec:processsetup}). Consequently, since the distance between the peaks of two species grows linearly with the column length while the width of each individual peak grows only with its square root, a longer column improves the separation of multiple species. Increasing the column length from 30\,cm to 120\,cm therefore yields an
increasingly clear separation of the six injected particle sizes, at the expense of broader \acp{rtd}.
Combining the inlet number based size distribution $\qCin(x)$ of the nanoparticles with the relationships shown in \cref{fig:fc} allows the full two-dimensional, size- and time-dependent outlet flux to be constructed. \Cref{fig:heatmap_qout} shows such a map for a column length of $60\;\mathrm{cm}$ and an injection duration of $1\;\mathrm{s}$. To investigate the chromatographic model independently of the synthesis stage, we prescribe a synthetic Gaussian inlet particle-number distribution with a mean diameter of $15\;\mathrm{nm}$, a standard deviation of $3\;\mathrm{nm}$, restricted to particle diameters between $1$ and
$23\;\mathrm{nm}$. Each horizontal line represents the temporal outlet signal associated with one particle size, weighted by the number of particles of this size present in the injection according to $\qCin(x)$. Conversely, each vertical line represents the mixture of particle sizes leaving the column at a given retention time. The inclined high-intensity region reflects the size-dependent retention behavior, i.e., smaller particles pass the column more slowly than their larger counterparts because of their higher pore accessibility.

\begin{figure}[pos=htbp]
    \centering
\includegraphics[width = \textwidth]{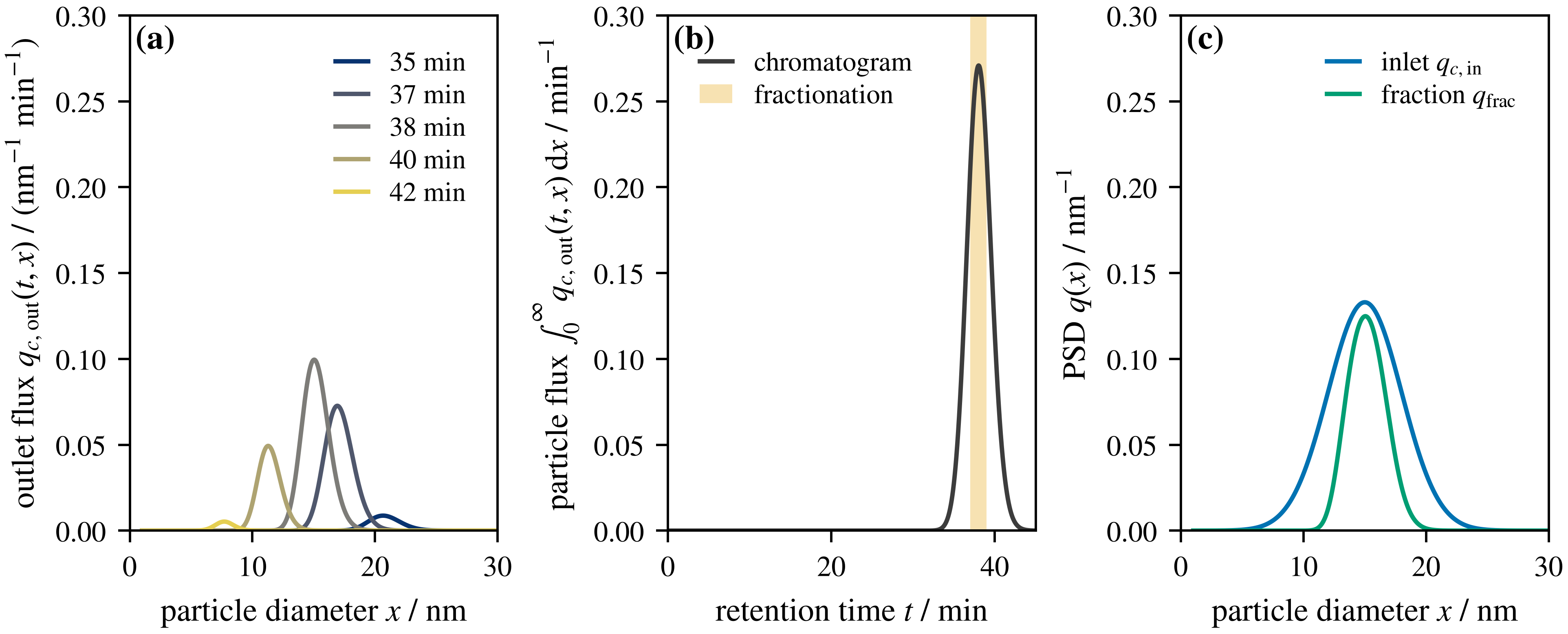}
    \caption{Fractionation mechanism for the reference chromatography case. (a) Size-dependent particle fluxes at retention times of $35$, $37$, $38$, $40$, and $42\,\mathrm{min}$. (b) Complete size-integrated chromatogram and the portion collected within $[t_0,t_1]=[37,39]\,\mathrm{min}$. (c) Absolute inlet and collected particle-number distributions.}
    \label{fig:qout_comparison}
\end{figure}

\Cref{fig:qout_comparison}a shows several vertical cuts through the two-dimensional distribution of \cref{fig:heatmap_qout}, i.e., the \acp{psd} exiting the column at retention times between 35 and 42\,min. Note that  these curves represent instantaneous size-resolved particle fluxes with units of $\mathrm{particles\,nm^{-1}\,min^{-1}}$. Consistent with \cref{fig:heatmap_qout}, high retention times such as 40\,min show a considerably lower number of particles, since the majority of the particles has already left the column by that time. As discussed in \cref{sec:processsetup,sec:separation}, size-exclusion chromatography can be used to narrow the \ac{psd} post-synthetically. A fractionation interval selects a subset of the particle population, typically excluding small and large particles from both edges of the \ac{psd} represented by the chromatogram. \Cref{fig:qout_comparison}b shows the complete size integrated chromatogram $\int_0 ^\infty \qCout(t,x) \mathrm d x$ in black. The orange area shows the portion collected between $t_0=37\,\mathrm{min}$ and $t_1=39\,\mathrm{min}$. Integrating the outlet flux over the fractionation interval (see \eqref{eq:qfrac} gives the collected \ac{psd} $\qfrac(x)$ shown in \cref{fig:qout_comparison}c. Note that although fractionation is intended to reduce the width of the \ac{psd}, it does so only indirectly, by effectively reducing the width of the chromatogram.

\subsection{Model for the coupled synthesis and separation process}
\label{sec:coupled_model}

\begin{figure}[pos=htbp]
      \centering
      \includegraphics[width=\textwidth]
      {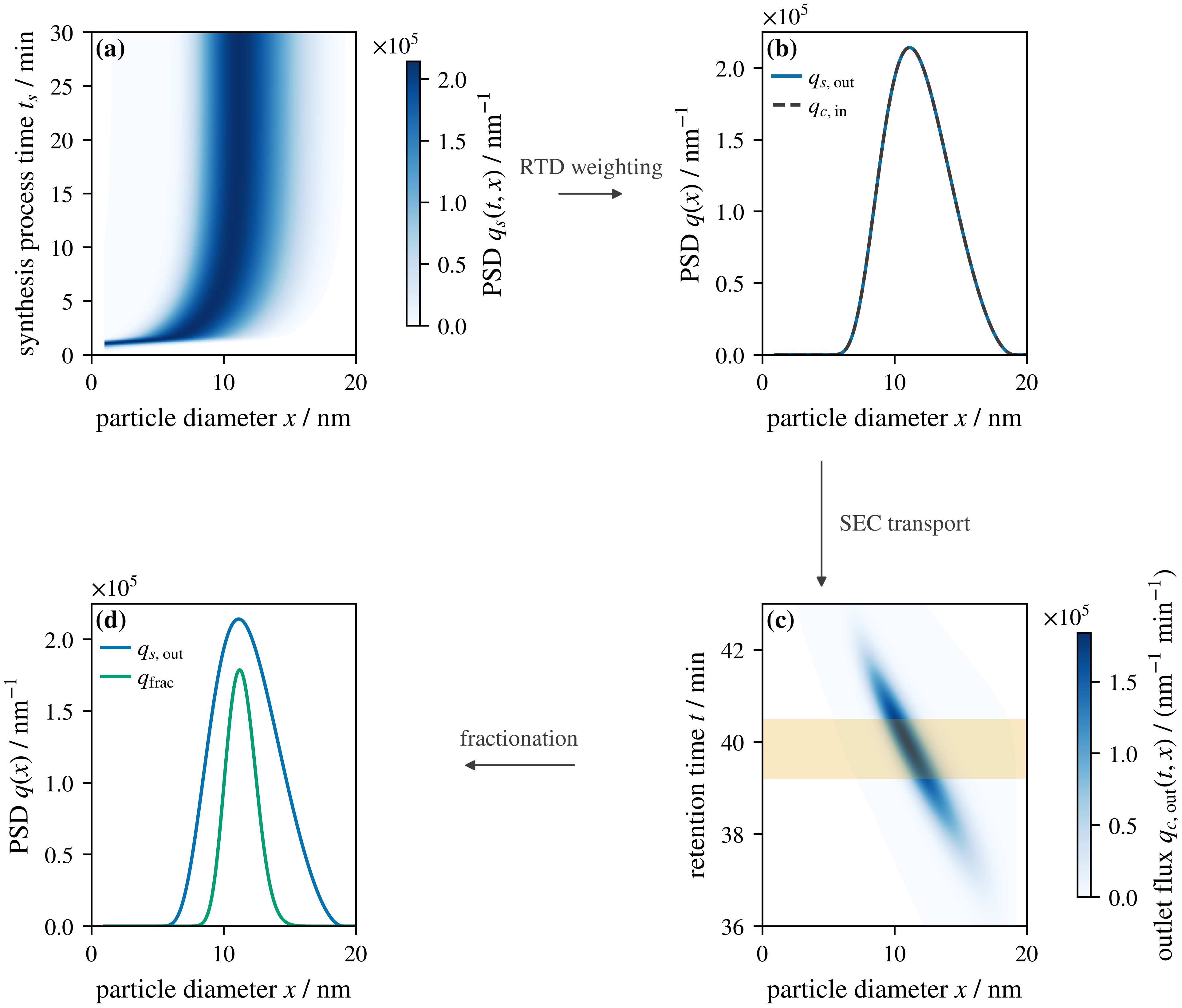}
      \caption{Evolution of the \ac{psd} through the coupled
      synthesis--separation process. (a) Time-dependent particle number
      distribution in one volume element during synthesis using $\kr = 0.2$. (b) Residence time weighted \ac{psd} at the reactor outlet with $\LS = 10 \mathrm m$, which is identical to the chromatographic inlet distribution,
      $\qSout(x)=\qCin(x)$. (c) Size- and time-dependent particle flux
      $\qCout(t,x)$ at the outlet of a chromatography column with $\LC=60\,\mathrm{cm}$ and injection duration of $1 \mathrm s$;
      the orange band delimits the fractionation interval
      $[t_0,t_1]=[39.2,40.5]\,\mathrm{min}$. (d) Reactor outlet \ac{psd} compared to the fractionated \ac{psd}s.}
      \label{fig:coupled_process_forward}
  \end{figure}

Having established the models for nanoparticle synthesis in a continuous flow reactor in \cref{sec:synthesis} and chromatographic separation by size-exclusion in
\cref{sec:separation}, we now combine both process steps into one process-chain model. As described in the process setup of \cref{sec:synthesis} we assume a continuous synthesis setup in which the particulate product is not prone to batch-to-batch variations. From the reactor outlet, the particles are transferred to a storage vessel and from there directly into the chromatographic column. Thus, for the following derivation, the reactor outlet \ac{psd} serves as the inlet \ac{psd} of
the chromatographic column, i.e.,
\[
    \qCin(x)=\qSout(x),
    \qquad x\geq\xn.
\]
Substituting this relation into the analytical fractionation model
\eqref{eq:qfracAnalytical} and further deriving $\qSout(x)$ as a solution of the synthesis model \eqref{eq:synthesismodel}, yields the coupled process model

\begin{align}
\label{eq:coupled_model}
\qfrac(x)
&=
\tfrac{\qSout(x)}{\tinj}
\Big[
    \HCcol(\tend,x;\LC)
    -\HCcol(\tstart,x;\LC) 
    -\HCcol(\tend-\tinj,x;\LC)
    +\HCcol(\tstart-\tinj,x;\LC)
\Big].
\end{align}
In \Cref{sec:synthesis} and \Cref{sec:separation} we have already established that our model reproduces the dynamics of both unit operations (see \Cref{fig:psd_k_L}-\ref{fig:qout_comparison}). When finally coupling both unit operations together and 
looking at the full process chain, it is important to understand that the \ac{psd} throughout the process chain constantly changes its shape and dimensionality due to the distribution of different particle sizes in time and space. \Cref{fig:coupled_process_forward} shows the most important steps in this dimensional evolution for a representative set of operating conditions. \newline

The process chain begins with the synthesis of nanoparticles in a continuous, tubular reactor. The synthesis is resolved on the level of a closed and perfectly mixed fluid element that travels through the tubular reactor and within which nucleation and growth proceed (see \Cref{sec:synthesis}). Such an element is characterized by its age, that is, by the time elapsed since it entered the reactor. Because nucleation and growth continuously redistribute mass over the size coordinate, the element does not possess a single \ac{psd} but a family of them, one for every age, forming a two-dimensional \ac{psd}, in particle size and process time (see \Cref{fig:coupled_process_forward}a). The reactor as a whole is then simulated by convoluting the single volume element with the \ac{rtd} of the tube, see \eqref{eq:qSout}. Under the assumption of complete segregation, the volume elements travel at the velocity dictated by their radial position so that the \acp{psd} differ based on the axial and radial position in the reactor. The dimensionality is therefore increased to a three-dimensional distribution. Observing the \ac{psd} at the reactor outlet specifically, fixes the axial coordinate at the reactor length, $z = \LS$, and records the mixture of all \acp{psd} along the diameter of the reactor. It thus leaves a single, one-dimensional \ac{psd} (see \Cref{fig:coupled_process_forward}b). The following intermediate storage tank downstream of the reactor does not alter this picture, since growth is assumed to cease upon leaving the reactor and the tank is treated as ideally mixed. \\

From the storage tank, the \ac{sec} column is fed by a pulsed injection and is consequently operated in a transient manner. Within the column, the particles are separated according to their access to the pore volume of the stationary phase, so that at a fixed point in time the \ac{psd} differs from position to position along the column. At the same time, and in contrast to the reactor, which is modeled in steady state, the \ac{psd} at a fixed position also changes with time, since particles are only present at that position while the pulse is passing. The description of the \ac{psd} within the column thus becomes three-dimensional again. \newline
Evaluating the \ac{psd} at the column outlet only, where $z = \LC$, again removes the axial coordinate of the column and leaves a two-dimensional \ac{psd} in particle size and retention time (see \Cref{fig:coupled_process_forward}c). In practice, however, what is recorded at the outlet is not this two-dimensional distribution but a one-dimensional chromatogram, i.e., a size-integrated signal (see, e.g., \cref{fig:qout_comparison}b) that is
typically weighted by the detector response, such as the extinction recorded by a diode array detector. At this point the dimensionality is thus further reduced. This step differs in kind from the two preceding reductions, because the dimensionality is not reduced by evaluating the distribution at a fixed
time, at a fixed spatial position, or at a single particle size, but by integrating over all particle sizes to obtain the chromatogram.

With knowledge of the system specific calibration curve, which connects the particles' diameter to their retention time, and the axial dispersion of the column, which blurs this relationship and must be deconvoluted, the chromatogram can finally be translated into a one-dimensional \ac{psd} of the classified product (see \Cref{fig:coupled_process_forward}d).\\

Thus, over the course of the process chain, the dimensionality of the data continuously changes as $2 \rightarrow 3 \rightarrow 1 \rightarrow 3 \rightarrow 2 \rightarrow 1 \rightarrow 1$-dimensional. 
  

%% file: 3_Optimization.tex

\section{Process-chain optimization}
\label{sec:optimization}
\subsection{Process design}
\label{sec:controls}
To formulate the process-structure relationship required for the subsequent optimization, we collect the process parameters in a unifying vector
\(
    \control\in\mathcal{V},
\) where $\mathcal{V}$ denotes the set of admissible process designs. Thus, the resulting \emph{process-structure mapping} can be summarized as
\begin{equation}
\label{eq:process2productmapping}
\begin{aligned}
    \ptop\colon\mathcal{V}
    &\longrightarrow \mathcal{Q},
    \\[-0.2em]
    \control
    &\longmapsto
    \ptop(\control)
    :=
    \qfrac(\,\cdot\,;\control).
    \\
\end{aligned}
\end{equation}

Here, $\mathcal{Q}$ denotes the set of resulting fractionated \acp{psd}.
Hence, $\ptop$ maps an admissible design of the complete process chain to the
corresponding fractionated \ac{psd} $\qfrac$ given by \eqref{eq:coupled_model}.
In this work, we consider two alternative process design formulations. In the first formulation, the reduction rate constant is optimized while the reactor length is fixed,\begin{equation*}
\label{eq:control_kr}
\control_{\kr}
=
\bigl(\kr,\LC,\tstart,\tend\bigr)^\top.
\end{equation*}
In the second formulation, the reactor length is optimized while the reduction rate constant is fixed,
\begin{equation*}
\label{eq:control_ls}
\control_{\LS}
=
\bigl(\LS,\LC,\tstart,\tend\bigr)^\top.
\end{equation*}
Thus, the two formulations differ only in the synthesis control, whereas the chromatographic column length and the fractionation interval are optimized in both cases. For either formulation $i\in\{\kr,\LS\}$, the corresponding admissible set is defined as
\begin{equation}
\label{eq:admissible_set}
\cV_i
:=
\left\{
\control_i\in\R^4
\,\middle|\,
\underline{\control}_{i,j}
\leq
\control_{i,j}
\leq
\overline{\control}_{i,j}
\quad \forall j,
\qquad
\tstart<\tend
\right\}.
\end{equation}
As demonstrated in \cref{sec:forward_sim_syn} and \cref{sec:forward_sim_chrom}, these process parameters strongly influence the resulting \ac{psd} and are thus relevant to be controlled. The proposed framework is, however, not restricted to
this particular set of parameters. Depending on the considered process, different or additional parameters may be included as tunable process parameters.

\subsection{Objective function}

As outlined in \cref{sec:intro}, nanoparticulate products are finding use in a growing number of technological applications. 
From the perspective of a nanoparticle manufacturer, a reasonable objective is to maximize overall economic profit $\profit$ which is the difference between the revenue $\revenue$ generated from selling it and the costs $\cost$ arising from its production and downstream processing:
\begin{equation*}
\label{eq:profit_general}
    \profit
    :=
    \revenue-\cost.
\end{equation*}

\paragraph{Revenue}
The revenue generated from selling the fractionated product depends on both the amount of recovered material and its quality. 
With a wide range of synthesis routes now established, the bottleneck is increasingly not access to nanoparticulate products per se, but rather to produce them at high quality. The ultimate, and still challenging, goal is to obtain nanoparticles with a narrow \ac{psd}, as broad distributions generally result in inferior product quality. For most fabrication processes, the goal is therefore to produce nanoparticles with small dispersity. Thus, product quality is characterized relative to a prescribed desired particle size~$\xd$. A high-quality product should therefore contain particles that are narrowly distributed around~$\xd$.
For a given set of process parameters
$\control\in\cV$, the process-to-structure mapping
\eqref{eq:process2productmapping} provides the fractionated \ac{psd} $\ptop(\control)$. We quantify its deviation from the desired particle size by the mean squared deviation
\begin{equation}
\label{eq:msdq}
    \msd^2(\xd;\control)
    :=
    \int_0^\infty
    (x-\xd)^2
    \frac{\ptop(\control)(x)}{\ym{0}(\control)}
    \,\D x,
\end{equation}
where
\begin{equation*}
    \ym{0}(\control)
    :=
    \int_0^\infty
    \ptop(\control)(x)\,\D x
\end{equation*}
denotes the zeroth moment of the final product distribution. 
Consequently, we obtain the 
relative root mean squared deviation by
\begin{equation*}
    \rmsd(\xd;\control) = \frac{\msd}{\xd}.
\end{equation*}

The proposed quality measure penalizes both the spread of the \ac{psd} and the deviations of its mean from the desired particle size. Consequently, if the mean particle size coincides with the desired size \(\xd\), the measure reduces to the standard deviation, and the corresponding relative measure coincides with the relative standard deviation.

Since $\qfrac(\cdot;\control)$ denotes the \ac{psd} collected within the fractionation interval, the corresponding recovered product mass is given by
\begin{equation*}
\label{eq:recovered_mass}
    M_{\mathrm{frac}}(\control)
    =
    \frac{\rho\pi}{6}
    \int_{\xn}^{\infty}
    x^3\ptop(\control)(x)\,\D x.
\end{equation*}
The overall process yield is defined as the recovered particle mass relative to the total precursor mass $p$ supplied to the synthesis process:
\begin{equation}
\label{eq:yield}
    \yield(\control)
    =
    \frac{M_{\mathrm{frac}}(\control)}{p}
    =
    \frac{\rho\pi}{6p}
    \int_{\xn}^{\infty}
    x^3\ptop(\control)(x)\,\D x.
\end{equation}
The yield can remain below unity even after precursor reduction is complete because monomer remaining in solution is not incorporated into the particles.
The notion of product quality is reflected in the price at which such products can be sold. A brief, non-exhaustive survey of gold nanoparticles on the market revealed that product price scales steeply with the relative standard deviation of the particles' core diameter. Based on our market survey (ESI, Figure 1), we assume an exponential relationship between the product's mean squared deviation from the target particle size and its price. (We emphasize that this concrete assumption could also be replaced by other functions with similar characteristics, without affecting our general approach.) The revenue of the coupled system is thus formulated as

\begin{equation}
\label{eq:revenue}
    \revenue(\control) = \eta_1 e ^{ - \eta_2  \rmsd(\xd;\control)} \yield(\control).
\end{equation}
Here, $\eta_1$ defines the revenue scale and $\eta_2>0$ determines the sensitivity to the dimensionless relative root mean square deviation. 
\paragraph{Costs}
For the given setup, costs can be broken down into precursor costs, costs for chromatography eluents, and  opportunity costs, which arise when the process takes longer than necessary and thereby delays restarting it. We note that there are also fixed costs, such as the initial investment in equipment; however, these act as constant offsets that do not affect the location of the optimum and can therefore be neglected in the optimization. \newline
Of the three contributions, precursor costs are accounted for indirectly through the yield \(\yield\), where a higher yield increases profit relative to the input costs, which discourages the algorithm from discarding most of the produced nanoparticles in pursuit of maximal purity. Furthermore, we consider costs only once the process has reached steady state, i.e., after reactor start-up. Opportunity costs for the synthesis therefore do not apply, since nanoparticles are produced continuously. Accordingly, we treat the chromatographic separation as the main contributor to process costs.

This contribution comprises two parts. The first is eluent consumption, i.e., the costs for the liquid mobile phase that transports the particle suspension through the column. The second is opportunity cost, as a longer separation cycle increases the time required per unit of product. For fixed operating conditions, the required eluent volume increases with the chromatographic process time. Since retention times scale linearly with column length (see \eqref{eq:mutr}), a longer column leads to correspondingly larger eluent consumption and reduces productivity at constant volumetric flow, because the separation takes longer.
All of these costs therefore scale with the column length \(\LC\) and can be collapsed into a single factor \( \costfac \), the cost per unit column length. 
A cost function for the process can then be formulated as

\begin{equation}
\label{eq:cost}
    \cost(\control; \costfac) = \costfac \LC.
\end{equation}

We emphasize again that the model is flexible: additional or different cost terms can be readily incorporated into the cost function. Treating the cost as linear in the column length is a reasonable first approximation, and can be replaced by a more detailed cost model if warranted. 

\paragraph{Overall optimization problem for coupled synthesis and separation}
\label{sec:overall_optimization}

The overall objective is to maximize the profit of the coupled process over
all admissible process parameters. Using the revenue and cost functions defined above,
the profit is given by
\begin{equation}
\label{eq:profit}
    \profit(\control;\costfac)
    :=
    \revenue(\control)-\cost(\control;\costfac).
\end{equation}
Here, the dependence of the revenue on the complete process chain is induced
by the process-to-product mapping $\ptop$.
Accordingly, the optimal process design is obtained by solving
\begin{equation}
\label{eq:overall_optimization}
    \max_{\control_i\in\cV_i}
    \profit(\control_i;\costfac),
    \qquad i\in\{\kr,\LS\},
\end{equation}
where the admissible sets $\cV_i$ are defined in \eqref{eq:admissible_set}.

The underlying process-design problem involves two competing criteria:
maximizing the revenue generated by product quality and yield, and
minimizing the cost of chromatographic separation. The profit function
\eqref{eq:profit} combines these criteria into a scalar objective and thus effectively forms the objective function for the subsequent optimization. By varying
the cost coefficient $\costfac$, different trade-offs between product
performance and separation cost can be explored.

In summary, in the following we optimize the process for a maximum in profit gained from selling the particulate product of the described process. We thus determine process parameters $\control_i$, with $i\in\{\kr,\LS\}$ such that the profit is maximized. 

\subsection{Inverse design of nanoparticle dispersity}
\label{sec:inversedesign}

Finally, we solve the optimization problem defined in
\eqref{eq:overall_optimization}. We investigate how the optimal process design responds to changes in the intrinsic performance of the two unit operations. In each case, the column cost coefficient $\costfac$ is varied
to trace different compromises between product revenue and separation cost.

The desired particle size is set to
\begin{equation*}
    \xd = 15~\mathrm{nm}.
\end{equation*}

\begin{table}[htbp]
\centering
\caption{Control variables and physical bounds used in the two optimization formulations.}
\label{tab:control_bounds}
\renewcommand{\arraystretch}{1.2}
{
\begin{tabular}{lcccc}
\hline
\textbf{Control} & \textbf{Symbol} & \textbf{Lower} & \textbf{Upper} & \textbf{SI-Unit}\\
\hline
Fractionation start    & $\tstart$ & $0$       & $9600$& s\\
Fractionation end      & $\tend$   & $0$       & $9600$& s\\
Column length & $\LC$   & $10^{-7}$ & $2$   & m\\
Reduction rate     & $\kr$     & $3.333 \times 10^{-4}$& $3.333 \times 10^{-3}$& s$^{-1}$\\
Reactor length    & $\LS$     & $2$       & $25$    & m\\
\hline
\end{tabular}}
\end{table}
Then the two cases introduced in section \ref{sec:forward_sim_syn} are examined again. In the first, $\kr$ is varied at constant $\LS$, and in the second, $\LS$ is varied at constant $\kr$. The bounds on $\kr$ apply only when optimizing with $\control_{\kr}$, whereas the bounds on $\LS$ apply only for $\control_{\LS}$. In both cases, the fractionation times additionally satisfy $\tstart<\tend$.

Furthermore, the parameters for the revenue function \eqref{eq:revenue} are given by \cref{tab:economic_parameters}, 
\begin{table}[t]
\centering
\caption{Economic model parameters.}
\label{tab:economic_parameters}
\begin{tabular}{llll}
\toprule
Parameter & Value  &Unit& Description\\
\midrule
$\eta_1$ & $363.33$  &€& Maximum normalized revenue\\
$\eta_2$ & $18.08$  &--& Exponential price sensitivity\\
$\gamma_\ell$ & $[0,\gamma_\ell^{\max}]$  &€ m$^{-1}$& Column-length cost coefficient\\
\bottomrule
\end{tabular}
\end{table}
whereas the upper bound for $\costfac$ can be approximately set by
\begin{equation}
\gamma_{\ell}^{\max}
=
\frac{\eta_1}{L_{\mathrm{c,max}}},
\end{equation}
where $L_{\mathrm{c,max}}$ is the upper bound of the chromatographic column length.

In this computational study we investigate how chromatographic efficiency changes
the role assigned to the separation stage in the coupled process. The size-dependent reference NTP value is scaled according to
\eqref{eq:ntp_scaling}, using $\alphaNTP\in\{0.25,1,2\}$ to represent different levels of chromatographic efficiency. Increasing
$\alphaNTP$ narrows the outlet signal of each particle size and results in a reduction of overlap between neighboring sizes. It determines how selectively the
fractionation interval can extract particles around $\xd$ from the supplied
$\qSout$.

\begin{figure}[pos=htbp]
    \centering
    \includegraphics[width=\textwidth]{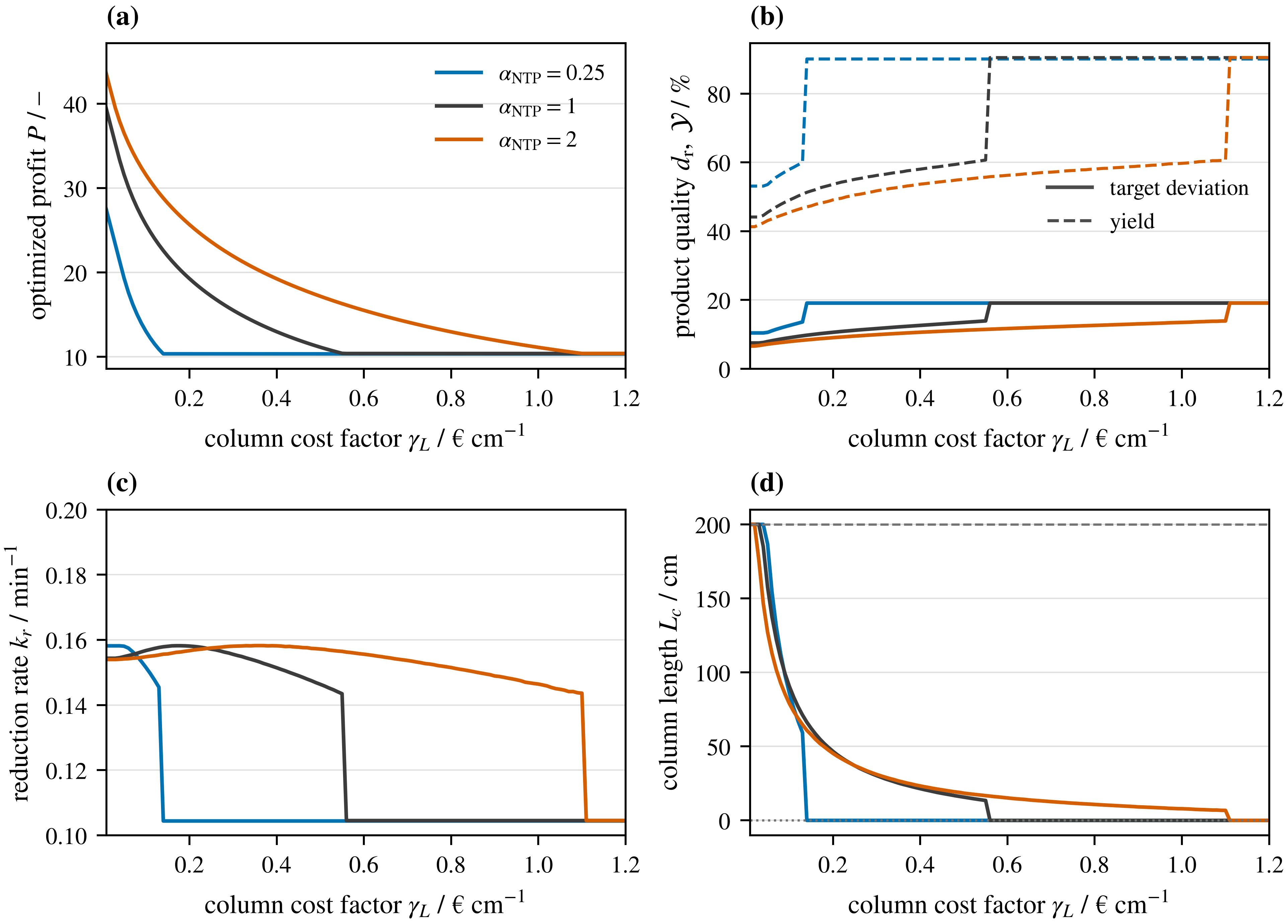}
    \caption{Locally optimized designs with $\kr$ as synthesis control for
    different separation efficiency scales $\alphaNTP$. Panels show
    (a) profit, (b) relative target deviation (solid) and yield (dashed),
    (c) reduction rate, and (d) chromatographic column length. Only
    $\costfac\geq10^{-2}\,\mathrm{cm^{-1}}$ is displayed.}
    \label{fig:optimization_kr_alpha_ntp}
\end{figure}

\begin{figure}[pos=htbp]
    \centering
    \includegraphics[width=\textwidth]{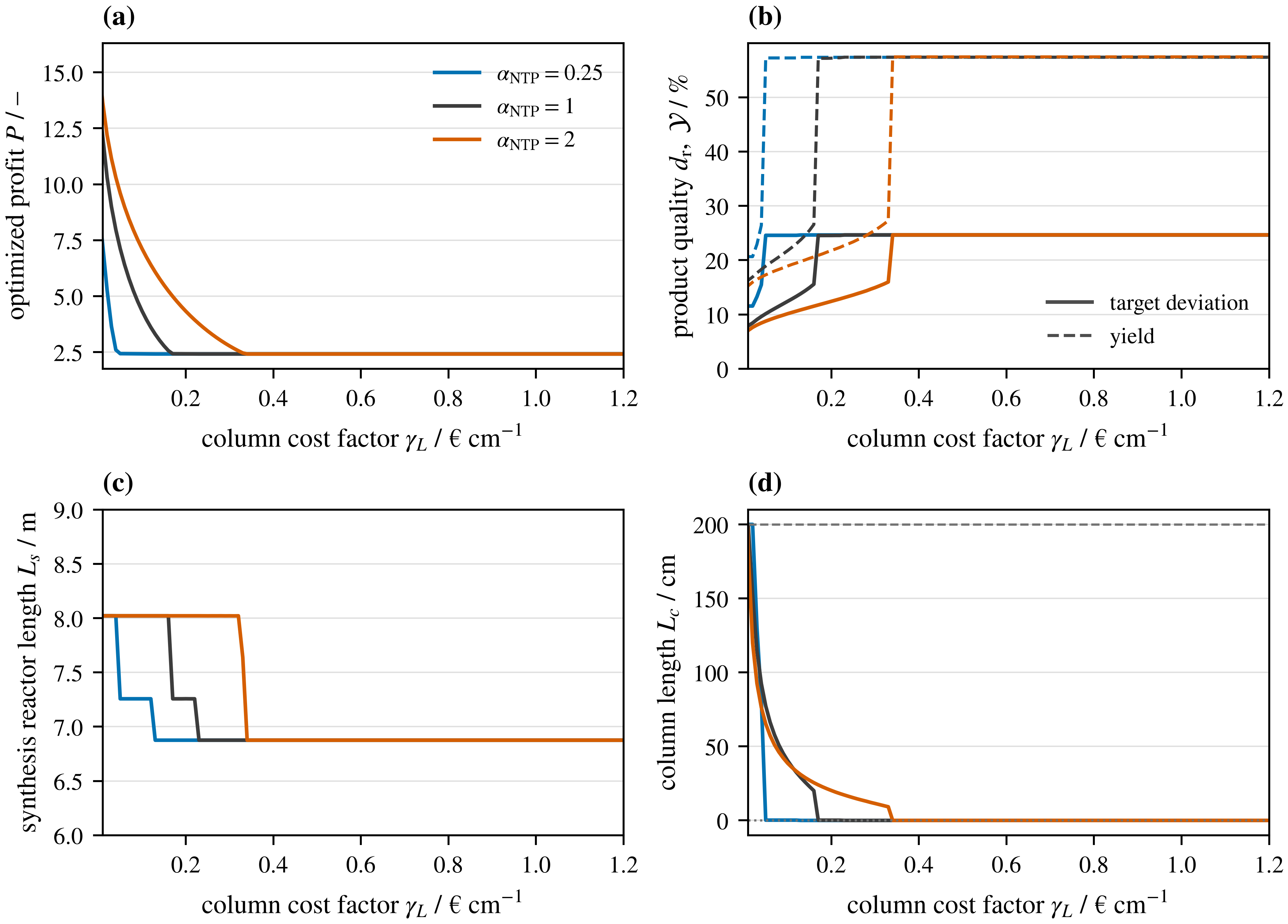}
    \caption{Locally optimized designs with $\LS$ as synthesis control for
    different separation efficiency scales $\alphaNTP$. Panels show
    (a) profit, (b) relative target deviation (solid) and yield (dashed),
    (c) synthesis reactor length, and (d) chromatographic column length. Only
    $\costfac\geq10^{-2}\,\mathrm{cm^{-1}}$ is displayed.}
    \label{fig:optimization_ls_alpha_ntp}
\end{figure}

For optimizing with $\control_{\kr}$ as well as with $\control_{\LS}$,
\Cref{fig:optimization_kr_alpha_ntp} and \Cref{fig:optimization_ls_alpha_ntp} show that (i) lower chromatographic costs favor longer columns and thus a higher chromatographic efficiency (see \Cref{fig:optimization_kr_alpha_ntp}d and \Cref{fig:optimization_ls_alpha_ntp}d) and (ii) higher chromatographic
efficiency increases product quality, thus also profit (see \Cref{fig:optimization_kr_alpha_ntp}a and \Cref{fig:optimization_ls_alpha_ntp}a), and influences the range over which chromatographic separation is applied. When $\alphaNTP$ is larger, fewer off-target particle sizes
elute together with the desired fraction. The optimizer can therefore choose a
more selective collection interval and obtain a smaller target deviation. As
already demonstrated by the collected distribution in
\cref{fig:qout_comparison}, this improvement is bought by collecting only part
of the inlet distribution. In the present optimization, improved
separation performance, impacted by separation efficiency due to $\alphaNTP$ as well as column length $\LC$, is deliberately used to favor quality over yield because the revenue increases strongly for a product concentrated around $\xd$.

\paragraph{Optimization with reduction rate constant as control for synthesis}

First, we analyze the results of the process-chain optimization using $\control_{\kr}$ with a fixed reactor length of $10\; \mathrm m$, see \Cref{fig:optimization_kr_alpha_ntp}. 
It provides information on the interplay of the unit operations of the process chain on different cost scales. 
If the costs for the chromatographic separation are low, a longer column provides sufficient separation to compensate for an outlet \ac{psd} whose mean particle size does not coincide with the desired size. The optimizer can therefore select a slightly higher reduction rate up to around $0.1582 \; \mathrm{min}^{-^1}$ (see \Cref{fig:optimization_kr_alpha_ntp}c), which speeds up the synthesis process, thus more mass can be produced at the fixed reactor length $\LS$. However, the mean of $\qSout$ is shifted away from $\xd$. The long column subsequently isolates a comparatively narrow fraction around the desired particle size. As $\costfac$ increases, the optimized column becomes shorter and its ability to correct the synthesis product decreases. The location of $\qSout$ then becomes increasingly important for the quality of the collected product. Consequently, the optimized reduction rate decreases down to $0.1434 \; \mathrm{min}^{-^1}$. This sacrifices part of the synthesis yield but shifts the mean particle size closer to $\xd$, thereby reducing the amount of chromatographic correction required. At sufficiently large chromatographic costs, the optimization ultimately favors a column-free process. In this regime, product quality must be achieved solely by the synthesis stage itself. At this point, the relative standard deviation and the yield of the process chain product increase to that of the synthesis product as no chromatographic separation is performed for column cost factors above this value (see \Cref{fig:optimization_kr_alpha_ntp}b).
 Representative distributions along this optimization path are shown in
  Figure~2 of the ESI. They illustrate how  $\qfrac$ approaches the synthesis outlet $\qSout$ as the optimized column becomes shorter. At the transition to the practically column-free design, the reduction in
  $\kr$ lowers the synthesis yield, but the collected yield increases because
  nearly the complete synthesis product is retained.

  \paragraph{Optimization with reactor length as control for synthesis}

In the second formulation, the reduction rate constant is fixed at $\kr = 0.05\,\mathrm{min}^{-1}$, while the synthesis reactor length $\LS$
is optimized together with the chromatographic controls. In contrast to $\kr$, varying $\LS$ does not alter the intrinsic nucleation and growth
kinetics, but determines how far the synthesis has progressed when the particles leave the reactor. Nevertheless, \Cref{fig:optimization_ls_alpha_ntp} reveals a similar, although not identical, interplay between the synthesis and separation controls.

As long as chromatographic separation is economically favorable, the optimizer selects a comparatively long reactor, thereby favoring a high synthesis yield despite a larger deviation of the reactor outlet \ac{psd} from the desired particle size $\xd$ (see \ref{fig:optimization_ls_alpha_ntp}c). With increasing column costs, the optimized reactor length decreases down to $7.26 \; \mathrm m$. This shifts the synthesis \ac{psd} towards the desired particle size while sacrificing part of the synthesis yield, thereby reducing the correction required from the chromatographic stage. The higher the chromatographic efficiency, the larger the column cost up to which this compensation by the separation stage remains economically feasible. Once chromatographic separation becomes too
expensive and the column is no longer used, all $\alphaNTP$ scenarios converge to the same reactor length at $6.88 \; \mathrm m$. In this regime, $\LS$ is selected such that the unseparated reactor outlet \ac{psd} provides the most favorable compromise between particle yield and proximity to the desired particle size. The intermediate levels of $\LS$ therefore represent changes between different process design strategies. Corresponding distributions for representative points of the $\LS$ optimization are provided in Figure~3 of the ESI. 
\newline

Across both control formulations, the dominant response to increasing column cost is a reduction in the optimized chromatographic column length
$\LC$. Higher chromatographic efficiency delays this transition because a selective product fraction remains economically favorable up to larger values of $\costfac$. Once chromatographic separation is abandoned, the influence of $\alphaNTP$ vanishes and all efficiency scenarios converge to the same synthesis-controlled design. The accompanying increase in yield and target deviation reflects the fundamental quality-yield trade off:
chromatographic fractionation improves product quality by rejecting off-target particles, but necessarily reduces the amount of collected
product.

\paragraph{Optimization procedure}
To handle the non-convex
optimization problem robustly and to obtain a smooth continuation of the optimized designs across the considered $\costfac$ values, the computational procedure was carried out in several successive refinement steps. First, promising starting points were identified by a coarse exploration of the admissible controls. These candidates were then refined by a local derivative-free coordinate search with decreasing step sizes. For
successive values of $\costfac$, the preceding solution was included as an additional starting point. For the $\LS$ control formulation, the result of the coordinate search was additionally refined using the sequential least-squares programming algorithm (SLSQP) implemented in SciPy. The reported designs are therefore locally optimized candidates, thus global optimality is not guaranteed.\newline
The computations were performed on an AMD Ryzen 7 PRO 8840U processor with eight physical cores and 16 threads, using Python 3.12, NumPy 2.5, and SciPy 1.18. Solving the synthesis model and reconstructing the corresponding reactor outlet \ac{psd} required approximately $1.2 \; \mathrm{s}$ per control setting on the specified hardware. Once the synthesis model is calculated, a coupled
forward evaluation required substantially less than one second, and the local refinement over the complete $\costfac$ path typically required less than one minute per scenario. The preceding coarse exploration accounted for the main computational effort.

%% file: 4_Conclusion.tex

\section{Conclusions and Outlook}
\label{sec:conclusion}

The computational experiments demonstrate the principal advantage of the coupled process-chain formulation. The quality of the final product is not determined by either synthesis or separation alone. 
Instead, the results reveal a transition from a separation dominated to a synthesis dominated process design. At low column costs, selective fractionation in a sufficiently long \ac{sec} column compensates for deviations of the synthesis outlet distribution from the desired particle size. As the column cost factor increases, the optimized column becomes shorter and a larger fraction of the synthesis product is collected. When chromatographic separation is no longer economically favorable, the collected \ac{psd} approaches the synthesis outlet \ac{psd}.

The position of this transition depends strongly on the chromatography separation efficiency. The results demonstrate that synthesis control, chromatographic performance, product quality, and yield cannot be considered independently. Their economically preferred balance changes with the relative cost of the separation step. The proposed framework addresses this challenge by evaluating both unit operations within a common optimization problem and explicitly accounting for their mutual influence on the final product.

The results emphasize the complexity of optimizing the process chain design and points out why sequential optimization of the individual unit operations can be suboptimal. Fixing the synthesis before designing the separation prevents the optimization from exploiting trade-offs across the process boundary. By contrast, the proposed coupled model directly selects the synthesis and separation controls with respect to the quality and economic performance of the final fractionated product.

Future research will focus on translating the approach developed in this work into laboratory practice. Furthermore, a subsequent step consists in developing an optimized, integrated process chain that is protected against uncertainty, as has already been performed separately for robust synthesis~\cite{Kuchlbauer2024} and robust separation~\cite{Cebulla}.